\documentclass{elsarticle}
\usepackage[final]{changes}
\newcommand{\stkout}[1]{\ifmmode\text{\sout{\ensuremath{#1}}}\else\sout{#1}\fi}
\setdeletedmarkup{\stkout{#1}}
\definechangesauthor[name={Roland Herzog},color=red]{RH}
\definechangesauthor[name={Martin Stoll},color=blue]{MS}
\definechangesauthor[name={John Pearson},color=violet]{JP}
\usepackage[utf8]{inputenc}
\usepackage{amsmath}
\usepackage{amsfonts}
\usepackage{graphicx}
\usepackage{mathrsfs}
\usepackage{mathenv}
\usepackage{amsopn}
\usepackage{amstext}
\usepackage{amssymb}
\usepackage{color}
\usepackage{subfig}
\usepackage{color}
\usepackage{bm}
\usepackage{multirow}
\usepackage[capitalize,nameinlink]{cleveref}
\usepackage{mathtools}
\usepackage[countmax]{subfloat}
\usepackage{pgfplots}
\pgfplotsset{compat=newest}
\pgfplotsset{plot coordinates/math parser=false}
\usepackage{url}
\usepackage{color}
\usetikzlibrary{plotmarks} 
\newlength\figureheight
\newlength\figurewidth 
\usetikzlibrary{plotmarks} 
\usetikzlibrary{external}
\usepgfplotslibrary{external} 
\newtheorem{theorem}{Theorem}
\newtheorem{lemma}[theorem]{Lemma}

\newproof{pf}{Proof}

\newcommand{\R}{\ensuremath{\mathbb{R}}}

\newcommand{\E}{\mathcal{E}}
\newcommand{\Am}{\mathcal{A}}
\newcommand{\A}{\mathcal{A}}
\newcommand{\Ll}{\mathcal{L}}

\newcommand{\M}{\mathcal{M}}

\newcommand{\J}{\mathcal{J}}

\newcommand{\G}{\mathcal{G}}
\newcommand{\I}{\mathbf{y}}
\newcommand{\m}{\mathbf{m}}
\newcommand{\aI}{\mathbf{p}}

\newcommand{\dOmega}{\textup{d}\Omega}
\newcommand{\dGamma}{\textup{d}\Gamma}
\newcommand{\dt}{\textup{d}t}
\renewcommand{\d}{\textup{d}}
\newcommand{\gmres}{\textsc{Gmres}}
\newcommand{\cgs}{\textsc{CGS}}

\begin{document}

\begin{frontmatter}

\title{Fast iterative solvers for an optimal transport problem}
% \tnotetext[mytitlenote]{Fully documented templates are available in the elsarticle package on \href{http://www.ctan.org/tex-archive/macros/latex/contrib/elsarticle}{CTAN}.}

%% Group authors per affiliation:
\author{Roland Herzog}
\address{Technische Universit\"{a}t Chemnitz, Faculty of Mathematics, Reichenhainer Straße~41, 09126 Chemnitz, Germany (\texttt{roland.herzog@mathematik.tu-chemnitz.de})}
\author{John W. Pearson}
\address{School of Mathematics, The University of Edinburgh, James Clerk Maxwell Building, The King's Buildings, Peter Guthrie Tait Road, EH9~3FD, United Kingdom (\texttt{j.pearson@ed.ac.uk})}
\author{Martin Stoll}
\address{Technische Universit\"{a}t Chemnitz, Faculty of Mathematics, Reichenhainer Straße~41, 09126 Chemnitz, Germany (\texttt{martin.stoll@mathematik.tu-chemnitz.de})}

\begin{abstract}
Optimal transport problems pose many challenges when considering their numerical treatment. We investigate the solution of a PDE-constrained optimisation problem subject to a particular transport equation arising from the modelling of image metamorphosis. We present the nonlinear optimisation problem, and discuss the discretisation and treatment of the nonlinearity via a Gauss--Newton scheme. We then derive preconditioners that can be used to solve the linear systems at the heart of the (Gauss--)Newton method. With the optical flow in mind, we further propose the reduction of dimensionality by choosing a radial basis function discretisation that uses the centres of superpixels as the collocation points. Again, we derive suitable preconditioners that can be used for this formulation.
\end{abstract}

\begin{keyword}
PDE-constrained optimisation, Saddle point systems, Time-dependent PDE-constrained optimisation, Preconditioning, Krylov subspace solver, Optical Flow, Optimal Transport
\end{keyword}

\end{frontmatter}

% \linenumbers
\section{Introduction}
The problem of optimal transport is a longstanding and active area of research in applied mathematics and the sciences \cite{villani2008optimal}. Its numerical treatment provides many challenges to the mathematical community (see \cite{kuzmin2010guide,Lev92} and the references therein).
Our goal in this paper is to discuss the solution of a PDE-constrained optimisation problem where the constraint is given by a transport equation.  In the field of PDE-constrained optimisation one typically wants to minimise an objective function with the constraints given by one or more PDEs \cite{book::IK08,book::FT2010}.

Much of our analysis for this formulation builds upon the previous work \cite{BHT10}, for which we wish to devise new iterative methods and discretisation schemes. The study of the original transport problem goes back to the 18th century but a modern formulation was given in \cite{kantorovich2006problem,kantorovitch1958translocation}. Recent developments include the seminal paper \cite{benamou2000computational}, where the problem is formulated as a fluid mechanics problem. A very similar formulation of minimising an objective function subject to a transport equation constraint is found in \textit{optical flow} (cf.\ for example \cite{asz,borzi2003optimal,bruhn2005lucas,haber2006multilevel,horn1981determining, mang2015inexact,mang2017lagrangian}), which models the apparent `motion' of an image.

In this manuscript, we examine the optimisation and discretisation of this problem, with a particular focus on the efficient solution of the linear systems that arise at the heart of the outer nonlinear solver. Our primary choice of nonlinear solver for the nonlinear optimisation problem is a Gauss--Newton scheme \cite{BHT10,HabAO00}, though we also consider methods based on a full Newton method. As one can typically follow the route of first performing the discretisation followed by deriving the appropriate optimality conditions, or vice versa, we discuss both approaches. We also briefly analyse a modified formulation of the classical transport model. Our main goal is the derivation of the linear system of equations followed by the introduction of suitable preconditioners that allow a parameter-robust solution of the linear system that is the computational heart of the nonlinear iteration. Such preconditioners have recently received much attention (cf. \cite{axelsson2016comparison,PSW11,PW11, Z10}). We then 
illustrate that the preconditioners proposed are efficient both for synthetic data as well as practical imaging examples.

One of the key bottlenecks when considering imaging application is the vast amount of data as the complexity is often prohibitively high due to the large number of pixels. A contribution of this paper is to consider replacing the sparse linear systems arising for the optical flow problem with smaller dense systems obtained when the system is discretised using radial basis function (RBF) techniques \cite{larsson2003numerical,Wendland}. For the centres of the RBFs locations we choose superpixel centres \cite{achanta2012slic}, which we assume do not change too drastically between an original (given) image and a target image that we wish to achieve through our PDE-constrained optimisation model. Strategies based on superpixels or supervoxels have recently been used to reduce the complexity of the methods, and we refer to \cite{amat2013fast,chang2013superpixel,donne2015fast} for details.

The paper is structured as follows. In \cref{sec::problem} we introduce the problem formulation considered in this work. \cref{sec::discr} introduces the discretisation of the optimisation problem and the constraint via a finite difference scheme. We discuss both discretise-then-optimise and optimise-then-discretise schemes. After introducing a modification to the problem formulation, we discuss two general preconditioning strategies in \cref{sec:Preconditioning}. We then introduce a discretisation using RBFs in \cref{sec::rbf}. In \cref{sec::Results} we present numerical experiments, for both finite difference and RBF discretisation, demonstrating the effectiveness of our discretisation and preconditioning approaches.
\section{The optimal transport problem}
\label{sec::problem}
The problem we examine in this paper is one of minimising the functional
\begin{align}\label{TranspOptical}
	\begin{split}
		\E(y,\vec{m})={}&\frac{1}{2\gamma}\int_{\Omega}(y(\vec{x},1)-y_1(\vec{x}))^2~\dOmega+\frac{\delta}{2}\int_{0}^{1}\int_{\Omega}(y(\vec{x},t)-\bar{y}(\vec{x},t))^2~\dOmega \, \dt \\
		&\quad\quad+\frac{\beta}{2}\int_{0}^{1}\int_{\Omega}(Q\vec{m}(\vec{x},t))^{2}~\dOmega \, \dt,
	\end{split}
\end{align}
where $\beta$ and $\gamma$ are (positive) parameters that can be understood as \emph{regularisation} or \emph{penalty parameters}. The parameter $\gamma$ is chosen in such a way that the computed state  $y(\vec{x},1)$ is close to the true final state $y_1$ at time $T=1$. Here, the velocity $\vec{m}$ represents a control variable, and $Q$ is a differential operator (possibly $Q=\text{blkdiag}(I,I)$ or $Q=\text{blkdiag}(\nabla,\nabla)$). The problem is solved on a space-time grid $(\vec{x},t):=([x_1,x_2],t)\in\Omega\times[0,1]$, where $\Omega \subset \R^2$ denotes the domain occupied by the image.
%(possibly $Q=\left[\begin{array}{cc}
%I & 0 \\ 0 & I \\
%\end{array}\right]$ or $Q=\left[\begin{array}{cc}
%\nabla & 0 \\ 0 & \nabla \\
%\end{array}\right]$). The problem is solved on a space-time grid $(\vec{x},t):=([x_1,x_2],t)\in\Omega\times[0,1]$.

\begin{figure}[htb!]
	\centering
%     \subfloat[Initial image $y_0$]{\includegraphics[width=0.45\textwidth]{matlab/00056}}\qquad
%     \subfloat[Target image $y_1$]{\includegraphics[width=0.45\textwidth]{matlab/00057}}
	\subfloat[Initial image $y_0$]{\includegraphics[width=0.25\textwidth]{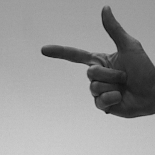}}\qquad
	\subfloat[Target image $y_1$]{\includegraphics[width=0.25\textwidth]{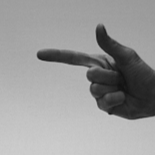}}
	\caption{An example of an optical flow problem where we have a starting picture on the left and a target picture on the right.}\label{fig::motivation}
\end{figure}

For the majority of the analysis presented in this paper we will consider the case $\delta=0$, i.e.\ where
\begin{align}\label{TranspCOST}
	\E(y,\vec{m})=\frac{1}{2\gamma}\int_{\Omega}(y(\vec{x},1)-y_1(\vec{x}))^2~\dOmega+\frac{\beta}{2}\int_{0}^{1}\int_{\Omega}(Q\vec{m}(\vec{x},t))^{2}~\dOmega \, \dt,
\end{align}
however on a number of occasions we will describe modifications which are taken into account when $\delta$ is a positive parameter, measuring the deviation of $y$ from the desired state $\bar{y}$ during the \textit{entire} time interval.
The goal is to minimise the above energy subject to the continuity transport equation
\begin{equation}
	\label{TranspEQN}
	y_t+\nabla\cdot (\vec{m}y) =0,
\end{equation}
with the initial condition $y(\vec{x},0)=y_0$ as well as appropriate boundary conditions, for instance periodic boundary conditions or Dirichlet conditions. Here, $\vec{m}=\left[m_1,~m_2\right]^T$ is defined for the two-dimensional domain $\Omega$.
While such a problem can be found in many areas of sciences, we wish to apply the above formulation to the estimation of an optical flow. To illustrate a particular set-up, examples for $y_0$ and $y_1$ are the two images shown in \cref{fig::motivation}\footnote{Images are taken from \url{http://cs.brown.edu/~black/images.html}.}.
\section{Discretisation using finite differences}
\label{sec::discr}
In this section, we wish to present how we discretise the optimisation problem \eqref{TranspCOST} with constraint \eqref{TranspEQN}. An outline is as follows: in \cref{sec::DTO} we examine the approach where such a PDE-constrained optimisation problem is discretised first, upon which optimality conditions are found. In \cref{sec::OTD} we then extend this methodology to the setting where optimality conditions are first derived (on a formal basis) on the continuous level, whereupon these are then discretised. In \cref{sec::for2} we then discuss the application of our methodology to a slight modification of the PDE \eqref{TranspEQN}. In \cref{sec:Modified}, we explain how the optimality conditions vary if one instead considers the cost functional \eqref{TranspOptical}, with an additional non-zero parameter $\delta$ measuring the deviation of the state $y$ from $\bar{y}$ throughout the entire time interval.

\subsection{Discretise-then-optimise}\label{sec::DTO}

A control problem using this formulation of the transport equation \eqref{TranspEQN} was introduced in \cite{BHT10}, and we therefore follow their approach for the derivation of the \textit{discretise-then-optimise} system. We start by discretising the objective function and nonlinear PDE constraint to then build a discrete Lagrangian, which then allows us to compute the solution via a Gauss-Newton or Lagrange--Newton scheme \cite[Ch.~10.3, 18]{nocedal1999no}. We employ an implicit Lax--Friedrichs method \cite{BHT10,Lev92} for the forward PDE
\begin{align*} 
	&\frac{1}{\tau}\left(y^{(k+1)}_{i,j}-\frac{1}{4}\left[y^{(k)}_{i+1,j}+y^{(k)}_{i-1,j}+y^{(k)}_{i,j+1}+y^{(k)}_{i,j-1}\right]\right)\\
	&\quad\quad+\frac{1}{2h}\left((m_1\odot y)^{(k+1)}_{i+1,j}-(m_1\odot y)^{(k+1)}_{i-1,j}+(m_2\odot y)^{(k+1)}_{i,j+1}-(m_2\odot y)^{(k+1)}_{i,j-1}\right)=0,
\end{align*}
that we can manipulate to arrive at the following system for each time-step:
\begin{equation*}
	\left(I+\frac{\tau}{2h}K(\m^{(k+1)})\right) \I^{(k+1)}=D_t\I^{(k)},\quad{}k=0,1,...,N_t-1.
\end{equation*}
Here, $\odot$ denotes the (componentwise) Hadamard product of vectors, $\tau$ represents the time-step and $h$ the spatial mesh parameter, the matrix $D_t$ arises from the four point stencil used to approximate the time derivative, and
\begin{equation*}
	K(\m^{(k+1)})=\left[\begin{array}{cc}
			D_1 & D_2 \\
	\end{array}\right]\left[
		\begin{array}{cc}
			\mathrm{diag}\big(\m_1^{(k+1)}\big)\\
			\mathrm{diag}\big(\m_2^{(k+1)}\big)\\
	\end{array}\right],
\end{equation*}
where $D_1$ and $D_2$ are centred finite difference matrices. We can then formulate an all-at-once approach using the notation 
\begin{equation*}
	L(\m^{(k+1)}) \, \I^{(k+1)}=D_t\I^{(k)},
\end{equation*}
where $L(\m^k)=I+\frac{\tau}{2h}K(\m^k)$, to obtain for a given number $N_t$ time-steps a matrix system of the form
\begin{equation*}
	\begin{array}{cccc}
		\underbrace{
			\left[
				\begin{array}{ccccc}
					L(\m^{(1)})&&&&\\
					-D_t&L(\m^{(2)})&&&\\
					&-D_t&\ddots &&\\
					&&\ddots &\ddots &\\
					&&&-D_t&L(\m^{(N_t)})\\
				\end{array}
		\right]}&
		\underbrace{
			\left[
				\begin{array}{c}
					\I^{(1)}\\
					\I^{(2)}\\
					\vdots\\
					\vdots\\
					\I^{(N_t)}\\
				\end{array}
		\right]}&=&
		\underbrace{
			\left[
				\begin{array}{c}
					D_t\I^{(0)}\\
					\mathbf{0}\\
					\vdots\\
					\vdots\\
					\mathbf{0}\\
				\end{array}
		\right]}\\
		\Am(\m)&\I&=&\mathbf{d}
	\end{array}
\end{equation*}
representing the discretised PDE constraint. 
Depending on the boundary conditions (Dirichlet or periodic), the matrices $L(\m^k)$ and $D_t$ need to be slightly modified in rows pertaining to boundary nodes.
In this work we apply periodic boundary conditions, in analogy to the work of Benzi, Haber and Taralli \cite{BHT10}.
Furthermore, we may approximate the objective function \eqref{TranspCOST} on the discrete level by
\begin{equation*}
	\bm{\mathcal{E}}(\I,\m)=\frac{1}{2\gamma}(\I^{(N_t)}-\I_1)^{T}M(\I^{(N_t)}-\I_1)+\frac{\beta \, \tau}{2}\mathbf{\m}^{T}\mathcal{M}\mathbf{\m},
\end{equation*}
where $\M=\mathrm{blkdiag}(W,\ldots,W)$, and $W$ is obtained from the discretisation of the term $\int_{\Omega}(Q\vec{m}(\vec{x},t))^{2} \, \dOmega$ (which could simply be a scaled identity operator). Note that, for simplicity, we have not included possible scalings of the individual $W$ matrices in $\M$ as these depend on the choice of discretisation performed in time. We now form the discrete Lagrangian for this problem 
\begin{equation*}
	\bm{\mathcal{L}}(\I,\m,\aI)=\bm{\mathcal{E}}(\I,\m)+\aI^{T}\mathcal{Q}\left(\A(\m) \, \I-\mathbf{d}\right),
\end{equation*}
where $\mathcal{Q}$ is a matrix allowing us to interpret the Lagrange multiplier $\aI$ as a grid function. For simplicity we assume that $\mathcal{Q}=\tau h^2 I,$ with $I$ the identity of the appropriate dimension.
Following \cite{BHT10}, the computation of the first order conditions
\begin{equation*}
	\bm{\mathcal{L}}_\I=\mathbf{0},\quad \bm{\mathcal{L}}_\m=\mathbf{0},\quad \bm{\mathcal{L}}_\aI=\mathbf{0},
\end{equation*}
leads to 
\begin{subequations}
	\label{forKKT}
	\begin{align}
		\label{forKKT1}\gamma^{-1}\M_{N_t}(\I-\I_{0,1})+\Am(\m)^{T}\mathcal{Q} \, \aI&=\mathbf{0},\\ 
		\label{forKKT2}\beta \, \tau \M \m+\J(\I)^{T}\mathcal{Q} \, \aI&=\mathbf{0},\\ 
		\label{forKKT3}\mathcal{Q} \left(\Am(\m) \, \I - \mathbf{d}\right)&=\mathbf{0},
	\end{align}
\end{subequations}
where $\I_{0,1}$ is a vector containing vectors of zeros for every time-step, apart from the final step which contains the vector $\I_1$. The matrix $\M_{N_t}$ contains zero blocks at every time-step, apart from the final time-step, which gives rise to an identy matrix scaled by $h^2$, denoted as $M$. Note that in optical flow applications one is often given an image for every time step, meaning the matrix $\M_{N_t}$ can be modified to one that does not contain zero diagonal blocks, and the vector $\I_{0,1}$ contains all time instances of these images. Further, $\J(\I)$ denotes the block diagonal matrix $\mathrm{blkdiag}\left(J(\I^{(1)}),\ldots,J(\I^{(N_t)})\right)$, where
\begin{equation*}
	\frac{\tau}{2h}\left[\begin{array}{cc}
			D_1\quad D_2 \\
	\end{array}\right]\left[
		\begin{array}{cc}
			\mathrm{diag}\left(\I^{(j)}\right)&0\\
			0&\mathrm{diag}\left(\I^{(j)}\right)\\
	\end{array}\right]=:J(\I^{(j)}),
\end{equation*}
at each time-step $j=1,2,...,N_t$.
The equations \eqref{forKKT} represents a nonlinear system, which we have to treat with a nonlinear optimisation scheme. We follow \cite{BHT10} and use a Gauss--Newton method for the solution of the first order conditions, which leads to the matrix system
\begin{align}\label{for1system1}
	\left[
		\begin{array}{ccc}
			\gamma^{-1}\M_{N_t}&0&\Am(\m)^{T}\mathcal{Q}\\
			0&\beta \, \tau \M &\J(\I)^{T}\mathcal{Q}\\
			\mathcal{Q}\Am(\m)&\mathcal{Q}\J(\I)&0\\
		\end{array}
	\right]
	\left[
		\begin{array}{c}
			\mathbf{s}_{\I}\\
			\mathbf{s}_{\m}\\
			\mathbf{s}_{\aI}
		\end{array}
	\right]=
	-\left[
		\begin{array}{c}
			\gamma^{-1}\M_{N_t}(\I-\I_{0,1})+\Am(\m)^{T}\mathcal{Q} \, \aI\\ 
			\beta \, \tau \M \m+\J(\I)^{T}\mathcal{Q} \, \aI\\ 
			\mathcal{Q} \left(\Am(\m) \, \I -\mathbf{d}\right)
		\end{array}
	\right]
\end{align}
at every step of the nonlinear iteration.

\subsection{Optimise-then-discretise}\label{sec::OTD}

We now highlight that it is also possible to follow the \textit{optimise-then-discretise} approach, where we commence by considering the continuous Lagrangian
\begin{equation*}
	\Ll(y,\vec{m},p)=\mathcal{E}(y,\vec{m})+\int_0^{1}\int_{\Omega} p \, (y_t+\nabla\cdot (\vec{m} y))~\dOmega \, \dt,
\end{equation*}
and then searching for the continuous first order conditions. Note that, for brevity, we have omitted the initial and boundary conditions within this \replaced[id=RH]{Lagrangian}{cost} functional; these are also accounted for and reappear in the optimality conditions. Proceeding formally, by considering the Fr\'{e}chet derivatives of the Lagrangian $\Ll$ with respect to $y$, $\vec{m}$ and $p$\added[id=RH]{,  and integration by parts}, we then obtain the conditions
\begin{subequations}
	\label{optKKT}
	\begin{align}
		\label{optKKT1}-p_t-\vec{m} \cdot \nabla{}p &=0,\\
		\label{optKKT2} \beta \, Q^{*}Q\vec{m}-y\nabla{}p&=\vec{0},\\
		\label{optKKT3}y_t+\nabla\cdot (\vec{m} y)&=0,
% \nonumber\frac{1}{\gamma} (\I(1,\cdot)-\I_1)&=-\aI(1)\\
	\end{align}
\end{subequations}
together with the initial condition for $y$, and the final-time condition
\begin{equation*}
	\frac{1}{\gamma} (y(\cdot,1)-y_1)=-p(1),
% \nonumber\int_{0}^{1}\int_{\partial\Omega}\aI (\left[\I d\m + \m d\I\right]\cdot \mathbf{n})&
\end{equation*}
corresponding to the adjoint equation. We have now established the continuous first order conditions for the optimal transport problem. Equations \eqref{optKKT} represent a nonlinear set of equations, which need to be augmented by boundary and initial conditions.  Equivalently, we can write this as $G(\vec{z})$  with $\vec{z}=\left[y,~m_1,~m_2,~p\right]^T$  and solve this nonlinear problem using a Gauss--Newton or Newton's method. The latter is the Lagrange--Newton or sequential quadratic programming (SQP) scheme and in each iteration we need to solve
\begin{equation*}
	G'(\vec{z}_k) \, \vec{s}_k=-G(\vec{z}_k)=:\vec{b},
\end{equation*}
where $\vec{s}_k=\vec{z}-\vec{z}_k$, until convergence of the method is achieved.
We now need to form the derivative of $G$ to solve the Newton problem, and obtain
\begin{subequations}
	\begin{align}
% \beta \, Q^{*}Qs_m-(\I\nabla s_{\aI})_x-(s_\I\nabla \aI)_x-(\I\nabla s_{\aI})_y-(s_\I\nabla \aI)_y &=0\\
		\label{Adjoint} -( s_{p})_{_t}-\vec{s}_{m}\cdot\nabla{}p-\vec{m}\cdot\nabla s_{p}  &=b_1,\\
		\label{Gradient} \beta \, Q^{*}Q\vec{s}_{m}-y\nabla s_{p}-s_y \nabla{}p&=\vec{b}_2,\\
		\label{Forward} (s_{y})_{_t}+\nabla\cdot (\vec{m} s_y)+\nabla\cdot (\vec{s}_{m}y)&=b_3.
	\end{align}
\end{subequations}
We examine the discretisation of this system of equations, starting with the treatment of the term
\begin{equation*}
	(s_{y})_{_t}+\nabla\cdot (\vec{m} s_y)
\end{equation*}
in the forward equation \eqref{Forward}, using the implicit Lax--Friedrichs scheme. This gives
\begin{align*} 
	&\frac{1}{\tau}\left(\big(s^{(k+1)}_{y}\big)_{i,j}-\frac{1}{4}\left[\big(s^{(k)}_{y}\big)_{i+1,j}+\big(s^{(k)}_{y}\big)_{i-1,j}+\big(s^{(k)}_{y}\big)_{i,j+1}+\big(s^{(k)}_{y}\big)_{i,j-1}\right]\right)\\
	&\quad\quad+\frac{1}{2h}\left((m_1\odot s_y)^{(k+1)}_{i+1,j}-(m_1\odot s_y)^{(k+1)}_{i-1,j}+(m_2\odot s_y)^{(k+1)}_{i,j+1}-(m_2\odot s_y)^{(k+1)}_{i,j-1}\right).
\end{align*}
Written in the same form as for the discretise-then-optimise approach, the matrices corresponding to the term at each time-step are $\left(I+\frac{\tau}{2h}K(\m^{(k+1)})\right) \mathbf{s}_\I^{(k+1)}-D_t\mathbf{s}_\I^{(k)}$, for $k=0,1,...,N_t-1$.
The discretisation of $\nabla\cdot (\vec{s}_{m}y)$ is performed analogously and we obtain
\begin{equation*}
	\nabla\cdot (\vec{s}_{m}y)=\frac{\partial(s_{m_1}y)}{\partial{}x_1}+\frac{\partial(s_{m_2}y)}{\partial{}x_2},
\end{equation*}
which leads to
\begin{equation*}
	\frac{1}{2h}\left((s_{m_1}\odot y)^{(k+1)}_{i+1,j}-(s_{m_1}\odot y)^{(k+1)}_{i-1,j}+(s_{m_2}\odot y)^{(k+1)}_{i,j+1}-(s_{m_2}\odot y)^{(k+1)}_{i,j-1}\right),
\end{equation*}
\textcolor{black}{and when taking into account the multiplication by the time-step $\tau$ results in 
	\begin{equation*}
		\frac{\tau}{2h}\left((s_{m_1}\odot y)^{(k+1)}_{i+1,j}-(s_{m_1}\odot y)^{(k+1)}_{i-1,j}+(s_{m_2}\odot y)^{(k+1)}_{i,j+1}-(s_{m_2}\odot y)^{(k+1)}_{i,j-1}\right).
\end{equation*}}
We write this in matrix form as $\bm{\mathcal{J}}(\I) \, \mathbf{s}_{\m}.$  
Consider now the discretisation of the terms arising from the continuous gradient equation \eqref{Gradient}. For the term $-s_y \nabla{}p=\big[-s_y \frac{\partial{}p}{\partial{}x_1},~-s_y \frac{\partial{}p}{\partial{}x_2}\big]^T$, we will obtain terms of the form
\begin{equation*}
%-s_y\nabla{}p=-s_y
%\left[
%\begin{array}{c}
%\frac{\partial{}p}{\partial{}x_1}\\
%\frac{\partial{}p}{\partial{}x_2}\\
%\end{array}
%\right]\quad\rightarrow\quad
%-(s_y)_{i,j}^{(k+1)}
	(s_y)_{i,j}^{(k+1)}\frac{1}{2h}\left[
		\begin{array}{c}
			\left(p_{i+1,j}^{(k+1)}-p_{i-1,j}^{(k+1)}\right)\\
			\left(p_{i,j+1}^{(k+1)}-p_{i,j-1}^{(k+1)}\right)\\
		\end{array}
	\right]
\end{equation*}
\textcolor{black}{at the $(k+1)$-th time-step, which can be written in matrix form (containing the time-step $\tau$) as
	\begin{equation*}
		\frac{\tau}{2h}\left[
			\begin{array}{c}
				\mathrm{diag}\left(D_1^{T}\aI^{(k+1)}\right)\\
				\mathrm{diag}\left(D_2^{T}\aI^{(k+1)}\right)
			\end{array}
		\right]
		\mathbf{s}_{\I}^{(k+1)}.
	\end{equation*}
}
In an analogous way, we may discretise the term $-y\nabla{}s_p=\big[-y \frac{\partial{}s_p}{\partial{}x_1},~-y \frac{\partial{}s_p}{\partial{}x_2}\big]^T$
%\begin{equation*}
%-y\nabla\cdot{}s_p=-y 
%\left[
%\begin{array}{c}
%(s_p)_x\\           
%(s_p)_y\\
%\end{array}
%\right]
%\end{equation*}
by
\begin{equation*}
	-y_{i,j}^{(k+1)}\frac{1}{2h}
	\left[
		\begin{array}{c}
			(s_p)_{i+1,j}^{(k+1)}-(s_p)_{i-1,j}^{(k+1)}\\
			(s_p)_{i,j+1}^{(k+1)}-(s_p)_{i,j-1}^{(k+1)}
		\end{array}
	\right],
\end{equation*}
which in block matrix form will lead to terms of the form \textcolor{black}{
	\begin{equation*}
		\frac{\tau}{2h}
		\left[
			\begin{array}{cc}
				\mathrm{diag}\left(\I^{(k+1)}\right)&0\\
				0&\mathrm{diag}\left(\I^{(k+1)}\right)\\
		\end{array}\right]
		\left[\begin{array}{c}
				D_1^{T}\\ D_2^{T} \\
		\end{array}\right]
		\mathbf{s}_{\aI}^{(k+1)},
% =:J(\I^{(k+1)}),
\end{equation*}}
abbreviated by $\J(\I)^{T}s_\aI$.
Finally, let us analyse the terms within the adjoint equation \eqref{Adjoint}.
% We here ignore the superscript ${(j)}$ as the block-structure of the discretisation will be obvious. 
The term
\begin{equation*}
	-\vec{s}_{m}\cdot\nabla{}p =-s_{m_1}\frac{\partial{}p}{\partial{x_1}}-s_{m_2}\frac{\partial{}p}{\partial{x_2}},
\end{equation*}
is approximated at time $t_{k+1}$ by
\begin{equation*}
	-(s_{m_1})^{(k+1)}_{i,j}\frac{1}{2h}\left(p_{i+1,j}^{(k+1)}-p_{i-1,j}^{(k+1)}\right)-(s_{m_2})^{(k+1)}_{i,j}\frac{1}{2h}\left(p_{i,j+1}^{(k+1)}-p_{i,j-1}^{(k+1)}\right),
\end{equation*}
or in matrix form
\begin{equation*}
	\left[
		\begin{array}{cc}
			\mathrm{diag}\left(D_1^{T}\aI^{(k+1)}\right)&\mathrm{diag}\left(D_2^{T}\aI^{(k+1)}\right)\\
	\end{array}\right] \mathbf{s}_\m^{(k+1)}=:\G(\mathbf{p}) \, \mathbf{s}_\m^{(k+1)}.
\end{equation*}
By now it is clear that the collection of all the previously discretised expressions results in a linear system similar to the matrix \eqref{for1system1}  obtained from the discretise-then-optimise, Gauss--Newton approach. The last ingredient needed is a discretised version of the discretised adjoint operator, i.e.,
\begin{equation*}
	-(s_p)_{_t}-\vec{m}\cdot\nabla{}s_p=-(s_p)_{_t}-m_1\frac{\partial{}s_p}{\partial{}x_1}-m_2\frac{\partial{}s_p}{\partial{}x_2}.
\end{equation*}
An implicit Lax--Friedrichs scheme again uses forward averaged differences in time and centred differences in space, leading to equations of the form
\begin{align*} 
	&-\frac{1}{\tau}\left((s_p)^{(k+1)}_{i,j}-\frac{1}{4}\left[(s_p)^{(k)}_{i+1,j}+(s_p)^{(k)}_{i-1,j}+(s_p)^{(k)}_{i,j+1}+(s_p)^{(k)}_{i,j-1}\right]\right)\\
	&\quad\quad-\frac{1}{2h}\left((m_{1})_{i,j}\left((s_p)^{(k+1)}_{i+1,j}-(s_p)^{(k+1)}_{i-1,j}\right)+(m_{2})_{i,j}\left((s_p)^{(k+1)}_{i,j+1}-(s_p)^{(k+1)}_{i,j-1}\right)\right),
\end{align*}
which in turn may be summarised by matrices $\left(I+\frac{\tau}{2h}L(\m^{(k+1)})\right) \mathbf{s}_\aI^{(k+1)}-D_t\mathbf{s}_\aI^{(k)}$. These may be assembled for all time-steps into a high-dimensional linear system $\mathcal{B}(\m)$ similar to $\A(\m)$. We have now discretised the PDE-constrained optimisation problem using the optimise-then-discretise approach. \textcolor{black}{We now have obtained a matrix system of the form
	\begin{equation}
		\label{system1}
		\left[
			\begin{array}{ccc}
				\gamma^{-1}\M_{N_t}&\G(\aI) &\mathcal{B}(\m) \\
				\G(\aI)^{T} &\beta \, \tau \M &\mathcal{J}(\I)^{T}\\
				\Am(\m)&\J(\I)&0\\
			\end{array}
		\right]
		\left[
			\begin{array}{c}
				\mathbf{s}_{\I}\\
				\mathbf{s}_{\m}\\
			\mathbf{s}_{\aI}\\\end{array}
		\right]
		=\mathbf{b}.
\end{equation}}
% Multiplying each of the optimality conditions by $\tau$, to achieve consistency with the matrix system \eqref{for1system1} obtained using the discretise-then-optimise method, we obtain the following linear system:
% \begin{equation}
% \label{system1}
% \left[
% \begin{array}{ccc}
% \gamma^{-1}\M_{N_t}&\tau\G(\aI) &\tau\mathcal{B}(\m) \\
% \tau\G(\aI)^{T} &\beta \, \tau \M &\tau\mathcal{J}(\I)^{T}\\
% \tau\Am(\m)&\tau\J(\I)&0\\
% \end{array}
% \right]
% \left[
% \begin{array}{c}
% \mathbf{s}_{\I}\\
% \mathbf{s}_{\m}\\
% \mathbf{s}_{\aI}\\\end{array}
% \right]
% =\mathbf{b}.
% \end{equation}
Note that we have not yet established that the discretisation of the adjoint equation above leads to the desired form that $\Am(\m)^{T}=\mathcal{B}(\m)$. For both matrices the diagonal blocks are of interest, and we will discuss these particular blocks now. For $\A(\m)^{T}$, we obtain for the crucial diagonal blocks that
\begin{align*}
	\left(K(\m^{(k+1)})\right)^{T}\mathbf{s}_{\aI}^{(k+1)}&=
	\left[
		\begin{array}{cc}
			\mathrm{diag}\left(\m_1^{(k+1)}\right)&
			\mathrm{diag}\left(\m_2^{(k+1)}\right)\\
	\end{array}\right]
	\left[
		\begin{array}{c}
			D_1^{T} \mathbf{s}_{\aI}^{(k+1)}\\
			D_2^{T} \mathbf{s}_{\aI}^{(k+1)}\\
	\end{array}\right]\\
	&=
	\mathrm{diag}\left(\m_1^{(k+1)}\right)D_1^{T} \mathbf{s}_{\aI}^{(k+1)}+
	\mathrm{diag}\left(\m_2^{(k+1)}\right)D_2^{T} \mathbf{s}_{\aI}^{(k+1)},
\end{align*}
whereupon applying $D_1^{T}=-D_1$ clearly leads to the desired form within $\mathcal{B}(\m)$. 
% { Is the final time condition represented properly? } {\color{blue}Martin: this might link in with my change above.}

We emphasise that the matrix system \eqref{system1} was obtained using a full Newton method, as opposed to the analysis for the discretise-then-optimise method for which the Gauss--Newton approach of \cite{BHT10} is applied. The main consequence of this change in the outer iteration is the appearance of the $\G(\aI)$ and $\G(\aI)^T$ blocks in \eqref{system1}. 
We also point out that in \eqref{for1system1} the scaling matrix $\mathcal{Q}$ was used for the Lagrange multiplier following \cite{BHT10}. Such a scaling could also be used to make system  \eqref{system1} ressemble the discretise-then-optimise approach more closely.

\subsection{Alternative problem formulation} \label{sec::for2}

Whereas we focus for the most part on the optimal transport problem given in \eqref{TranspCOST}--\eqref{TranspEQN}, we also wish to briefly discuss an alternative formulation given by the minimisation of
\begin{align*}
	\E(y,\vec{m})=\frac{1}{2\gamma}\int_{\Omega}(y(\vec{x},1)-y_1)^2~\dOmega+\frac{\beta}{2}\int_{0}^{1}\int_{\Omega}(Q\vec{m}(\vec{x},t))^{2}~\dOmega \, \dt
\end{align*}
subject to the advection transport equation (cf. \cite{mang2017lagrangian})
\begin{equation}
	\label{TranspEQN2}
	y_t+\vec{m}\cdot \nabla{}y =0,
\end{equation}
along with appropriate boundary and initial conditions.
Let us briefly compare \eqref{TranspEQN2} with \eqref{TranspEQN}. The divergence theorem implies $\frac{\d}{\dt} \int_\Omega y \, \dOmega = - \int_{\Gamma} (\vec{m} y) \cdot \vec{n} \, \dGamma$ and thus \eqref{TranspEQN} will be mass preserving \deleted[id=RH]{up} in the presence of homogeneous Dirichlet or periodic boundary conditions. By constrast, mass may be produced or removed in \eqref{TranspEQN2} unless $\nabla\cdot\vec{m} = 0$ holds.

Discretising the objective function as before results in the following functional on the discrete level:
\begin{equation}\label{Objective_New}
	\bm{\mathcal{E}}(\I,\m)=\frac{1}{2\gamma}(\I^{(N_t)}-\I_1)^{T}M(\I^{(N_t)}-\I_1)+\frac{\beta \, \tau}{2}\mathbf{\m}^{T}\mathcal{M}\mathbf{\m}.
\end{equation}
The discretisation of the transport equation \eqref{TranspEQN2} via an implicit Lax--Friedrichs scheme \cite{Lev92} gives
\begin{align*} 
	&\frac{1}{\tau}\left(y^{(k+1)}_{i,j}-\frac{1}{4}\left[y^{(k)}_{i+1,j}+y^{(k)}_{i-1,j}+y^{(k)}_{i,j+1}+y^{(k)}_{i,j-1}\right]\right)\\
	&\quad\quad+\frac{1}{2h}\left((m_1)_{i,j}\left(y^{(k+1)}_{i+1,j}-y^{(k+1)}_{i-1,j}\right)+(m_2)_{i,j}\left(y^{(k+1)}_{i,j+1}-y^{(k+1)}_{i,j-1}\right)\right)=0,
\end{align*}
which can be written in matrix form as 
\begin{align*}
	\left(I+\frac{\tau}{2h}\widetilde{K}(\m^{(k+1)})\right) \I^{(k+1)}=D_t\I^{(k)},
\end{align*}
where 
\begin{equation*}
	\widetilde{K}(\m^{(k+1)})=
	\left[
		\begin{array}{cc}
			\mathrm{diag}\left(\m_1^{(k+1)}\right)&\mathrm{diag}\left(\m_2^{(k+1)}\right)\\
	\end{array}\right]
	\left[
		\begin{array}{c}
			D_1\\
			D_2\\
	\end{array}\right].
\end{equation*}
Therefore, in block form, the system of equations for the forward problem at all time-steps reads
\begin{equation}\label{Constraints_New}
	\begin{array}{cccc}
		\underbrace{
			\left[
				\begin{array}{ccccc}
					\widetilde{L}(\m^{(1)})&&&&\\
					-D_t&\widetilde{L}(\m^{(2)})&&&\\
					&-D_t&\ddots &&\\
					&&\ddots &\ddots &\\
					&&&-D_t&\widetilde{L}(\m^{(N_t)})\\
				\end{array}
		\right]}&
		\underbrace{
			\left[
				\begin{array}{c}
					\I^{(1)}\\
					\I^{(2)}\\
					\vdots\\
					\vdots\\
					\I^{(N_t)}\\
				\end{array}
		\right]}&=&
		\underbrace{
			\left[
				\begin{array}{c}
					M\I^{(0)}\\
					0\\
					\vdots\\
					\vdots\\
					0\\
				\end{array}
		\right]}\\
		\widetilde{\Am}(\m)&\I&=&\mathbf{d}
	\end{array}.
\end{equation}
Using the standard Lagrangian approach for differentiating the objective functional \eqref{Objective_New} subject to the constraints \eqref{Constraints_New}, we obtain the first order conditions
\begin{align*}
	\bm{\mathcal{L}}_\I&=\gamma^{-1}\mathcal{M}_{N_t}(\mathbf{y}-\I_{0,1})+\widetilde{\Am}(\m)^{T}\aI=\mathbf{0},\\
	\bm{\mathcal{L}}_\m&=\beta \, \tau \M \m+\widetilde{\G}(\I) \, \aI=\mathbf{0},\\
	\bm{\mathcal{L}}_\aI&=\widetilde{\Am}(\m) \, \I - \mathbf{d}=\mathbf{0},
\end{align*}
where
\begin{equation*}
	\widetilde{G}^{(j)}=
	\left[
		\begin{array}{c}
			\mathrm{diag}(\I_1^{(j)})\\
			\mathrm{diag}(\I_2^{(j)})
		\end{array}
	\right]~~\textnormal{with}~~\I_i^{(j)}=D_iy^{(j)}~~\textnormal{and}~~\widetilde{\G}=\mathrm{blkdiag}\left(\widetilde{G}^{(1)},\ldots,\widetilde{G}^{(N_t)}\right).
\end{equation*}
As for the previous problem formulation, we may then write down a Gauss--Newton scheme for this problem governed by the matrix
\begin{equation}
	\label{for2system1}
	\left[
		\begin{array}{ccc}
			\gamma^{-1}\M_{N_t}&0&\Am(\m)^{T}\\
			0&\beta \, \tau \M &\widetilde{\G}(\I)\\
			\mathcal{A}(\m)&\widetilde{\G}(\I)^T&0\\
		\end{array}
	\right].
\end{equation}

\subsection{Modified cost functional}\label{sec:Modified}

We now briefly discuss the changes to the optimality conditions and matrix systems if the cost functional \eqref{TranspOptical} is instead considered. In this case, when the discretise-then-optimise method is applied, the discrete approximation of $\mathcal{E}$ is given by
\begin{equation*}
	\bm{\mathcal{E}}(\I,\m)=\frac{1}{2\gamma}(\I^{(N_t)}-\I_1)^{T}M(\I^{(N_t)}-\I_1)+\frac{\delta\tau}{2}(\I-\bar{\I})^{T}\bar{\M}(\I-\bar{\I})+\frac{\beta \, \tau}{2}\mathbf{\m}^{T}\mathcal{M}\mathbf{\m},
\end{equation*}
where $\bar{\I}$ contains the discrete values of the desired state $\bar{y}$ at each time-step, and $\bar{\M}$ is a block diagonal matrix corresponding to a scaled identity operator applied at each time-step. The equations $\bm{\mathcal{L}}_{\mathbf{m}}=\bm{\mathcal{L}}_{\mathbf{p}}=\mathbf{0}$, as given by \eqref{forKKT2}--\eqref{forKKT3}, will then hold as before. By contrast, the equation $\bm{\mathcal{L}}_{\mathbf{y}}=\mathbf{0}$ becomes
\begin{equation*}
	\ \gamma^{-1}\M_{N_t}(\I-\I_{0,1})+\delta\tau\bar{\M}(\I-\bar{\I})+\mathcal{A}(\m)^{T}\mathcal{Q} \, \aI=\mathbf{0},
\end{equation*}
and the Gauss--Newton system \eqref{for1system1} is therefore modified to the form
\begin{align*}
	&\left[
		\begin{array}{ccc}
			\gamma^{-1}\M_{N_t}+\delta\tau\bar{\M}&0&\Am(\m)^{T}\mathcal{Q}\\
			0&\beta \, \tau \M &\J(\I)^{T}\mathcal{Q}\\
			\mathcal{Q}\Am(\m)&\mathcal{Q}\J(\I)&0\\
		\end{array}
	\right]
	\left[
		\begin{array}{c}
			\mathbf{s}_{\I}\\
			\mathbf{s}_{\m}\\
			\mathbf{s}_{\aI}
		\end{array}
	\right] \\
	&\quad\quad\quad\quad=
	-\left[
		\begin{array}{c}
			\gamma^{-1}\M_{N_t}(\I-\I_{0,1})+\delta\tau\bar{\M}(\I-\bar{\I})+\Am(\m)^{T}\mathcal{Q} \, \aI\\ 
			\beta \, \tau \M \m+\J(\I)^{T}\mathcal{Q} \, \aI\\ 
			\mathcal{Q} \left(\Am(\m) \, \I -\mathbf{d}\right)
		\end{array}
	\right].
\end{align*}

Similarly, for the matrix system \eqref{system1} arising from the optimise-then-discretise approach, the top left entry $\gamma^{-1}\M_{N_t}$ must be replaced by $\gamma^{-1}\M_{N_t}+\delta\tau\bar{\M}$ if the modified cost functional \eqref{TranspOptical} is used.

\section{Preconditioning}\label{sec:Preconditioning}
The most important step within our algorithm, in order to minimise the computational work required, is to accurately and efficiently solve large and sparse linear systems. To illustrate our methodology, we focus our description on the Gauss--Newton matrix derived from the discretise-then-optimise case for problem \eqref{TranspCOST}--\eqref{TranspEQN}; see \eqref{for1system1}:
\begin{equation}\label{Matrix}
	\left[
		\begin{array}{ccc}
			\gamma^{-1}\M_{N_t}&0&\Am(\m)^{T}\mathcal{Q}\\
			0&\beta \, \tau \M &\J(\I)^{T}\mathcal{Q}\\
			\mathcal{Q}\Am(\m)&\mathcal{Q}\J(\I)&0\\
		\end{array}
	\right].
\end{equation}
We approach the solutions of such linear systems by exploiting the \emph{saddle point form} of the matrices involved. It is well known that non-singular saddle point matrices of the form
\begin{equation}\label{SaddlePt}
	\ \left[\begin{array}{cc}
			A & B^T \\
			C & 0 \\
	\end{array}\right]
\end{equation}
can be effectively approximated (provided that $A$ is non-singular) by block diagonal or block triangular preconditioners of the form
\begin{equation*}
	\ P_D=\left[\begin{array}{cc}
			A & 0 \\
			0 & S \\
	\end{array}\right],\quad\quad{}P_T=\left[\begin{array}{cc}
			A & 0 \\
			B & -S \\
	\end{array}\right],
\end{equation*}
where $S:=CA^{-1}B^T$ denotes the (negative) \emph{Schur complement} of the matrix system. It can be shown (see \cite{Ku95,preconMGW}) that preconditioning the saddle point system with $P_D$ or $P_T$ results in the convergence of a Krylov subspace method in $3$ or $2$ iterations, respectively. It can also be shown (see \cite{ipsen}) that a similar block triangular preconditioner may also be applied if the $(2,2)$ block of \eqref{SaddlePt} is non-zero.

Of course, the so-called `ideal preconditioners' $P_D$ and $P_T$ would not be applied in practice, as the computational cost of inverting $A$ and $S$ would be almost as great as inverting the entire system. We therefore wish to consider variants of these preconditioners, where the $(1,1)$ block and Schur complement are replaced with suitable (cheap) approximations.
For the matrix \eqref{Matrix}, we see that
\begin{equation*}
	\ A=\left[\begin{array}{cc}
			\gamma^{-1}\M_{N_t}&0\\
			0&\beta \, \tau \M\\
	\end{array}\right],\quad\quad{}B=C=\left[\begin{array}{cc}
			\mathcal{Q}\Am(\m)&\mathcal{Q}\J(\I)\\
	\end{array}\right].
\end{equation*}
Below, we present two preconditioners which we discover to be effective for this system. In the case of `full observations', i.e., $\delta > 0$ and with the desired state given for every time-step, as in some optical flow problems, the matrix $\M_{N_t}$ is block diagonal and invertible, which makes preconditioning easier. We therefore focus on the case where $\M_{N_t}$ is highly singular and only comment on the more straightforward case.

\subsection{First preconditioner}
\label{sec::precond1}

The first preconditioner we introduce is based on the block diagonal structure $P_D$, but with the $(1,1)$ block and Schur complement replaced with suitable approximations. For the $(1,1)$ block, we write that
\begin{equation*}
	A=\left[\begin{array}{ccc}
			\gamma^{-1}\M_{N_t}&0\\
			0&\beta \, \tau \M \\
	\end{array}\right]\approx
	\left[\begin{array}{ccc}
			\widehat{\M}_{N_t}&0\\
			0&\beta \, \tau \M \\
	\end{array}\right]=:\widehat{A},
\end{equation*}
where $\widehat{\M}_{N_t}$ approximates the highly singular matrix $\gamma^{-1}\M_{N_t}$. As suggested by Benzi, Haber and Taralli in \cite{BHT10}, all zero diagonal entries in the $(1,1)$ block are replaced with $\mu$ within the preconditioner, where this parameter reflects the mean ratio of diagonal entries between the first and second terms of the Schur complement.

Since the $(1,1)$ block is highly singular we define the Schur complement of the ``perturbed'' system as
\begin{equation*}
	S=\mathcal{Q}\Am(\m)\widehat{\M}_{N_t}^{-1}\Am(\m)^{T}\mathcal{Q}+\frac{1}{\beta \, \tau}\hspace{0.15em}\mathcal{Q}\J(\I) \M^{-1}\J(\I)^{T}\mathcal{Q}.
\end{equation*}
The approach we use to approximate this matrix follows the matching strategy introduced in \cite{PS12_CP,PSW11,PW10,PW11}, where we approximate $S$ by 
\begin{equation}\label{S1}
	\widehat{S}_1=\big(\mathcal{Q}\Am(\m)+\mathcal{M}_1\big)\widehat{\M}_{N_t}^{-1}\big(\Am(\m)^{T}\mathcal{Q}+\mathcal{M}_2\big),
\end{equation}
with the desire that the term $\mathcal{M}_1\widehat{\M}_{N_t}^{-1}\mathcal{M}_2$ accurately captures the second term of the exact Schur complement, that is:
\begin{equation*}
	\mathcal{M}_1\widehat{\M}_{N_t}^{-1}\mathcal{M}_2\approx\frac{1}{\beta \, \tau}\hspace{0.15em}\mathcal{Q}\J(\I)\M^{-1}\J(\I)^{T}\mathcal{Q}.
\end{equation*}
A possible selection of the matrices $\mathcal{M}_1$, $\mathcal{M}_2$ is as follows:
\begin{equation}
	\ \label{M1M2} \mathcal{M}_1=\mathcal{M}_2^T=\frac{1}{\sqrt{\beta \, \tau}}\hspace{0.15em}\big(\mathcal{Q}\J(\I) \M^{-1}\J(\I)^{T}\mathcal{Q}\big)^{1/2}\widehat{\M}_{N_t}^{1/2},
\end{equation}
A further saving in the required computational cost may be achieved by replacing these matrices by the diagonal approximations:
\begin{equation*}
	\ \mathcal{M}_1=\mathcal{M}_2=\frac{1}{\sqrt{\beta \, \tau}}\hspace{0.15em}\left[\text{diag}\big(\widehat{\M}_{N_t}\big)\right]^{1/2}\left[\text{diag}\big(\mathcal{Q}\J(\I)\M^{-1}\J(\I)^{T}\mathcal{Q}\big)\right]^{1/2},
\end{equation*}
which leads to a preconditioner that may be applied cheaply in practice. We note that, whereas the matrices $\widehat{\M}_{N_t}^{1/2}$ and $\M^{-1}$ appear complicated to apply, they are in fact straightforward as each of the matrices $\M_{N_t}$ and $\M$ contain multiples of identity operators on each diagonal block, and therefore $\widehat{\M}_{N_t}^{1/2}$ and $\M^{-1}$ solely consist of scaled identity matrices corresponding to each time-step.

Whereas the effectiveness of this Schur complement approximation will inevitably depend to some extent on the numerical behaviour of the solution at each Newton step, the following observation may be readily made (based on the methodology of \cite{PW11}), guaranteeing the robustness of the smallest eigenvalue of the preconditioned Schur complement in an ideal setting:
\begin{lemma}
	The eigenvalues $\lambda$ of $\widehat{S}_1^{-1}S$ satisfy $\lambda\geq\frac{1}{2},$ \textcolor{black}{ where $\widehat{S}_1$ is as defined by \eqref{S1}, and $\mathcal{M}_1$, $\mathcal{M}_2$ are given in \eqref{M1M2}.}
\end{lemma}
\begin{pf}
 Due to the symmetry and positive definiteness of $S$ and $\widehat{S}_1,$ which may be observed due to $\widehat{\M}_{N_t}$ being symmetric positive definite by construction, the eigenvalues may be bounded by the Rayleigh quotient
\begin{equation*}
	\ R:=\frac{\mathbf{v}^T S \mathbf{v}}{\mathbf{v}^T\widehat{S}_1\mathbf{v}}=\frac{\mathbf{v}^T\left(\mathcal{Q}\Am(\m)\widehat{\M}_{N_t}^{-1}\Am(\m)^{T}\mathcal{Q}+\beta^{-1}\tau^{-1}\mathcal{Q}\J(\I) \M^{-1}\J(\I)^{T}\mathcal{Q}\right)\mathbf{v}}{\mathbf{v}^T\big(\mathcal{Q}\Am(\m)+\mathcal{M}_1\big)\widehat{\M}_{N_t}^{-1}\big(\Am(\m)^{T}\mathcal{Q}+\mathcal{M}_2\big)\mathbf{v}}.
\end{equation*}
% where $\mathcal{L}=\beta^{-1/2}\tau^{-1/2}\mathcal{Q}\J(\I)\M^{-1/2}\widehat{\M}_{N_t}^{1/2}$.
We now observe that we may write
\begin{equation*}
	\ R=\frac{\mathbf{a}^T\mathbf{a}+\mathbf{b}^T\mathbf{b}}{(\mathbf{a}+\mathbf{b})^T(\mathbf{a}+\mathbf{b})},~~~\text{where}~~\mathbf{a}=\widehat{\M}_{N_t}^{-1/2}\Am(\m)^{T}\mathcal{Q}\mathbf{v},~\mathbf{b}=\frac{1}{\sqrt{\beta \, \tau}}\hspace{0.15em}\M^{-1/2}\J(\I)^{T}\mathcal{Q}\mathbf{v}.
\end{equation*}
Simple manipulation therefore tells us that
\begin{equation*}
	\ \frac{1}{2}(\mathbf{a}-\mathbf{b})^T(\mathbf{a}-\mathbf{b})\geq0~~\Leftrightarrow~~\mathbf{a}^T\mathbf{a}+\mathbf{b}^T\mathbf{b}\geq\frac{1}{2}(\mathbf{a}+\mathbf{b})^T(\mathbf{a}+\mathbf{b})~~\Leftrightarrow~~R\geq\frac{1}{2},
\end{equation*}
which leads to the result.
\end{pf}

\vspace{0.5em}

\begin{figure}[htb!]
 \setlength\figureheight{0.24\linewidth} 
	\setlength\figurewidth{0.36\linewidth}
	\centering
	\subfloat[$n_x=8$]{
		\includegraphics[width=0.45\textwidth]{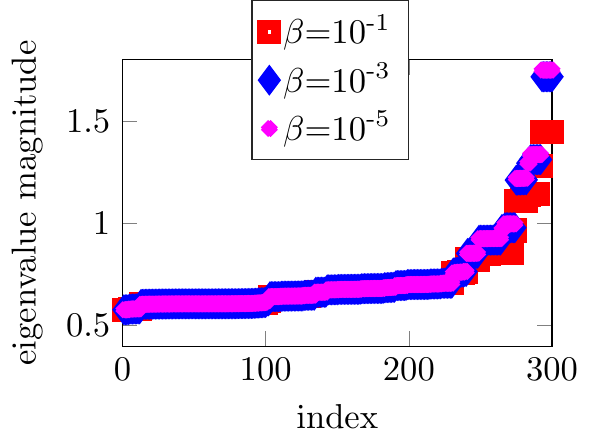}
	}
	\subfloat[$n_x=16$]{
		\includegraphics[width=0.45\textwidth]{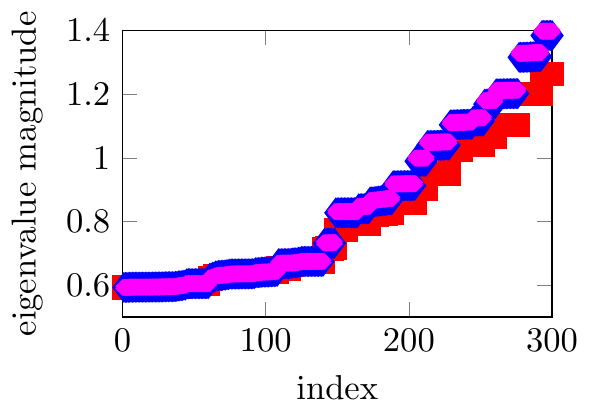}
	}
	\caption{Eigenvalues of $\widehat{S}_1^{-1}S$ for an imaging test problem, where finite difference nodes are equally distributed in each spatial direction. Results are given for different numbers of nodes $n_x$ in each dimension, and for different values of $\beta$.}\label{fig::eig}
\end{figure}

We highlight that, although the lower bound for the eigenvalues of $\widehat{S}_1^{-1}S$ can be analysed in detail, the magnitude of the largest eigenvalue will depend on the precise behaviour of the Newton iterates, which we cannot control in general. To provide an illustration of the overall eigenvalue distribution, we present in \cref{fig::eig} the eigenvalues for a particular test problem, for a range of problem sizes and values of $\beta$. As is demonstrated by the plots, the eigenvalues are found to become more clustered for finer grids, with the magnitude of the largest eigenvalues fairly robust to changes in regularisation parameter.

Applying our approximations of the $(1,1)$ block and Schur complement leads to a preconditioner of the form
\begin{equation}
	\label{eq:P1}
	\ P_1=\left[\begin{array}{ccc}
			\widehat{\M}_{N_t}&0&0\\
			0&\beta \, \tau\M&0\\
			0&0&\widehat{S}_1\\
	\end{array}\right],
\end{equation}
which we can then apply within a Krylov subspace method.

\vspace{0.5em}

We highlight that similar ideas may be applied to the matrix system \eqref{system1} arising from the optimise-then-discretise setting. In more detail, one may apply preconditioners of the form
\begin{equation*}
	\ \left[\begin{array}{ccc}
			\gamma^{-1}\M_{N_t}&0 &0 \\
			0 &\beta \, \tau \M &0\\
			0&0&\widehat{S}_{1,\text{OTD}}\\
	\end{array}\right]\quad\text{or}\quad\left[\begin{array}{ccc}
			\gamma^{-1}\M_{N_t}&0 &0 \\
			\G(\aI)^{T} &\beta \, \tau \M &0\\
			\Am(\m)&\J(\I)&-\widehat{S}_{1,\text{OTD}}\\
	\end{array}\right].
\end{equation*}
Here, $\widehat{S}_{1,\text{OTD}}$ can be chosen to approximate the Schur complement of the matrix system obtained by setting $\G(\I)=0$, for example. That is,
\begin{align*}
	\ \widetilde{S}:={}&\Am(\m)\widehat{\M}_{N_t}^{-1}\mathcal{B}(\m)+\frac{1}{\beta \, \tau}\J(\I)\M^{-1}\J(\I)^T \\
	\approx{}&\left(\Am(\m)+\frac{1}{\sqrt{\beta \, \tau}}\M_{1,\text{OTD}}\right)\widehat{\M}_{N_t}^{-1}\left(\mathcal{B}(\m)+\frac{1}{\sqrt{\beta \, \tau}}\M_{2,\text{OTD}}\right)=:\widehat{S}_{1,\text{OTD}},
\end{align*}
with $\M_{1,\text{OTD}}$, $\M_{2,\text{OTD}}$ chosen such that $\M_{1,\text{OTD}}\widehat{\M}_{N_t}^{-1}\M_{2,\text{OTD}}\approx\J(\I)\M^{-1}\J(\I)^T$.

\vspace{0.5em}

\textbf{Note.}~The derivation for this preconditioner has been based on the cost functional \eqref{TranspCOST}, whereupon the highly singular matrix $\gamma^{-1}\M_{N_t}$ is approximated by $\widehat{\M}_{N_t}$. If instead the cost functional \eqref{TranspOptical} is considered (with $\delta>0$), the corresponding matrix in the $(1,1)$ block of \eqref{Matrix} is $\M_{1,1}:=\gamma^{-1}\M_{N_t}+\delta\tau\bar{\M}$, which is now invertible. Therefore, when deriving an analogous preconditioner $P_1$ for this problem setup, the matrix $\widehat{\M}_{N_t}$ must be replaced with $\M_{1,1}$ on all occasions.

\subsection{Second preconditioner}
\label{sec::precond2}

We now derive a second block preconditioner, based largely on results in \cite{BenDOS15}. We commence by considering the following permutation of the matrix to be solved:
\begin{equation}\label{Permuted}
	\Pi
	\left[
		\begin{array}{ccc}
			\gamma^{-1}\M_{N_t}&0&\Am(\m)^{T}\mathcal{Q}\\
			0&\beta \, \tau \M &\J(\I)^{T}\mathcal{Q}\\
			\mathcal{Q}\Am(\m)&\mathcal{Q}\J(\I)&0\\
		\end{array}
	\right]=
	\left[
		\begin{array}{ccc}
			\mathcal{Q}\Am(\m)&\mathcal{Q}\J(\I)&0\\
			0&\beta \, \tau \M &\J(\I)^{T}\mathcal{Q}\\
			\gamma^{-1}\M_{N_t}&0&\Am(\m)^{T}\mathcal{Q}\\
		\end{array}
	\right],
\end{equation}
where the permutation matrix is given by
\begin{equation*}
	\Pi:=
	\left[
		\begin{array}{ccc}
			0&0&I\\
			0&I&0\\
			I&0&0\\
		\end{array}
	\right].
\end{equation*}
The matrix \eqref{Permuted} is now of saddle point structure \eqref{SaddlePt}, with
\begin{equation*}
	\ A=\left[
		\begin{array}{cc}
			\mathcal{Q}\Am(\m)&\mathcal{Q}\J(\I)\\
			0&\beta \, \tau\M\\
		\end{array}
	\right],\quad\quad{}B=\left[
		\begin{array}{cc}
			0&\mathcal{Q}\J(\I)\\
		\end{array}
	\right],\quad\quad{}C=\left[
		\begin{array}{cc}
			\gamma^{-1}\M_{N_t}&0\\
		\end{array}
	\right],
\end{equation*}
and a non-zero $(2,2)$ block given by $\Am(\m)^{T}\mathcal{Q}$.

We may then consider the right preconditioner
\begin{equation*}
	\widetilde{P}=
	\left[
		\begin{array}{ccc}
			\mathcal{Q}\Am(\m)&\mathcal{Q}\J(\I)&0\\
			0&\beta \, \tau \M &0\\
			\gamma^{-1}\M_{N_t}&0&-\widehat{S}_2\\
		\end{array}
	\right],
\end{equation*}
with its inverse given by
\begin{equation*}
	\widetilde{P}^{-1}=
	\left[
		\begin{array}{ccc}
			\Am(\m)^{-1}\mathcal{Q}^{-1}&-\beta^{-1}\tau^{-1}\Am(\m)^{-1}\J(\I) \M^{-1}&0\\
			0&\beta^{-1}\tau^{-1} \M^{-1} &0\\
			\gamma^{-1}\widehat{S}_2^{-1}\M_{N_t}\Am(\m)^{-1}\mathcal{Q}^{-1}&-\gamma^{-1}\beta^{-1}\tau^{-1}\widehat{S}_2^{-1}\M_{N_t}\Am(\m)^{-1}\J(\I) \M^{-1}&-\widehat{S}_2^{-1}\\
		\end{array}
	\right].
\end{equation*}
The matrix $\widehat{S}_2$ is designed to approximate the Schur complement of the \emph{permuted matrix system}, that is
\begin{equation*}
	\widehat{S}_2\approx
	S
% \left[
% \begin{array}{cc}
% \gamma^{-1}\M_{N_t}&0\end{array}
% \right]
% \left[
% \begin{array}{cc}
% \Am(\m)^{-1}\mathcal{Q}^{-1}&-\beta^{-1}\tau^{-1}\Am(\m)^{-1}\J(\I)\M^{-1}\\
% 0&\beta^{-1}\tau^{-1} \M^{-1}\\
% \end{array}
% \right]
% \left[
% \begin{array}{c}
% 0\\
% \J(\I)^{T}\mathcal{Q}\\
% \end{array}
% \right]\\
	=\Am(\m)^{T}\mathcal{Q}+\frac{1}{\beta \, \tau\gamma}\hspace{0.15em}\M_{N_t}\Am(\m)^{-1}\J(\I)\M^{-1}\J(\I)^{T}\mathcal{Q}.
\end{equation*}
%If we now shift the permutation $\Pi$ onto the preconditioner we obtain
%\begin{equation*}
%\widetilde{P}^{-1}\Pi=
%\left[
%\begin{array}{ccc}
%0&-\beta^{-1}\tau^{-1}\Am(\m)^{-1}\mathcal{Q}^{-1}\mathcal{Q}\J(\I) \M^{-1}&\Am(\m)^{-1}\mathcal{Q}^{-1}\\
%0&\beta^{-1}\tau^{-1} \M^{-1} &0\\
%-\widehat{S}^{-1}&0&\gamma^{-1}\widehat{S}^{-1}\M_{N_t}\Am(\m)^{-1}\mathcal{Q}^{-1}\\
%\end{array}
%\right]
%\end{equation*}
%and permuted back we obtain
Let us now reapply the permutation to the preconditioned system (that is to say we propose a preconditioner $P_2$ such that $P_2^{-1}=\widetilde{P}^{-1}\Pi$), and therefore obtain
\begin{equation}
	\label{eq:prec1}
	P_2^{-1}=
	\left[
		\begin{array}{ccc}
			0&-\beta^{-1}\tau^{-1}\Am(\m)^{-1}\J(\I)\M^{-1}&\Am(\m)^{-1}\mathcal{Q}^{-1}\\
			0&\beta^{-1}\tau^{-1} \M^{-1}&0\\
			-\widehat{S}_2^{-1}&-\gamma^{-1}\beta^{-1}\tau^{-1}\widehat{S}^{-1}\M_{N_t}\Am(\m)^{-1}\J(\I) \M^{-1}&\gamma^{-1}\widehat{S}_2^{-1}\M_{N_t}\Am(\m)^{-1}\mathcal{Q}^{-1}\\
		\end{array}
	\right].
\end{equation}
Applying the preconditioner is in fact more straightforward than it currently appears. To compute a vector $\mathbf{v}=P_2^{-1}\mathbf{w}$, where $\mathbf{v}:=\left[\mathbf{v}_{1}^T,~\mathbf{v}_{2}^T,~\mathbf{v}_{3}^T\right]^T$, $\mathbf{w}:=\left[\mathbf{w}_{1}^T,~\mathbf{w}_{2}^T,~\mathbf{w}_{3}^T\right]^T$, we first see from the second block of $P_2^{-1}$ that
\begin{equation*}
	\beta^{-1}\tau^{-1} \M^{-1}\mathbf{w}_2=\mathbf{v}_2.
\end{equation*}
The first equation derived from \eqref{eq:prec1} then gives that
\begin{align*}
	-\beta^{-1}\tau^{-1}\Am(\m)^{-1}\J(\I)\M^{-1}\mathbf{w}_2+\Am(\m)^{-1}\mathcal{Q}^{-1}\mathbf{w}_3={}&\mathbf{v}_1\\
	\Rightarrow~~-\Am(\m)^{-1}\J(\I)\mathbf{v}_2+\Am(\m)^{-1}\mathcal{Q}^{-1}\mathbf{w}_3=\Am(\m)^{-1}\big(\mathcal{Q}^{-1}\mathbf{w}_3-\J(\I)\mathbf{v}_2\big)={}&\mathbf{v}_1,
\end{align*}
and using this we can write the last equation in \eqref{eq:prec1} as
\begin{align*}
% -\widehat{S}^{-1}w_1-\gamma^{-1}\widehat{S}^{-1}\M_{N_t}\beta^{-1}\tau^{-1}\Am(\m)^{-1}\mathcal{Q}^{-1}\mathcal{Q}\J(\I) \M^{-1}w_2+\gamma^{-1}\widehat{S}^{-1}\M_{N_t}\Am(\m)^{-1}\mathcal{Q}^{-1}w_3={}&v_3\\
	&-\widehat{S}_2^{-1}\mathbf{w}_1-\gamma^{-1}\widehat{S}_2^{-1}\M_{N_t}\big(\beta^{-1}\tau^{-1}\Am(\m)^{-1}\J(\I) \M^{-1}\mathbf{w}_2-\Am(\m)^{-1}\mathcal{Q}^{-1}\mathbf{w}_3\big)=\mathbf{v}_3\\
	\Rightarrow~~&-\widehat{S}_2^{-1}\mathbf{w}_1+\gamma^{-1}\widehat{S}_2^{-1}\M_{N_t}\big(\Am(\m)^{-1}\mathcal{Q}^{-1}\mathbf{w}_3-\Am(\m)^{-1}\J(\I)\mathbf{v}_2\big) \\
	&\hspace{21.3em}=\widehat{S}_2^{-1}\big(\gamma^{-1}\M_{N_t}\mathbf{v}_1-\mathbf{w}_1\big)=\mathbf{v}_3.
\end{align*}
Therefore, in order to solve a system with the preconditioner $P_2$, we need to solve for the matrix $\M$, which is certainly invertible, as well as the matrix $\mathcal{Q}\Am(\m).$ What remains is the construction of the approximation $\widehat{S}_2$ of the Schur complement. In more detail, we suggest the use of a similar matching strategy as above, to write
\begin{align*}
% \Am(\m)^{T}\mathcal{Q}+\beta^{-1}\tau^{-1}\gamma^{-1}\M_{N_t}\Am(\m)^{-1}\mathcal{Q}^{-1}\mathcal{Q}\J(\I)\M^{-1}\J(\I)^{T}\mathcal{Q}
	\widehat{S}_2=\big(\Am(\m)^{T}\mathcal{Q}+\mathcal{M}_l\big)\Am(\m)^{-1}\mathcal{Q}^{-1}\big(\mathcal{Q}\Am(\m)+\mathcal{M}_r\big),
\end{align*}
where 
\begin{equation*}
	\mathcal{M}_l\Am(\m)^{-1}\mathcal{Q}^{-1}\mathcal{M}_r\approx\frac{1}{\gamma\beta \, \tau}\hspace{0.15em}\M_{N_t}\Am(\m)^{-1}\J(\I)\M^{-1}\J(\I)^{T}\mathcal{Q}.
\end{equation*}
Such an approximation may be achieved if, for example, 
\begin{align*}
	\mathcal{M}_l&= \frac{1}{\gamma}\hspace{0.15em}\M_{N_t},\\
	\mathcal{M}_r&= \frac{1}{\beta \, \tau}\hspace{0.15em}\mathcal{Q}\J(\I)\M^{-1}\J(\I)^{T}\mathcal{Q}\quad\text{or}\quad\frac{1}{\beta \, \tau}\hspace{0.15em}\mathcal{Q}\J(\I)\hspace{0.1em}\text{diag}(\M)^{-1}\J(\I)^{T}\mathcal{Q},
\end{align*}
and we thus build such approximations into our preconditioner $P_2$. We highlight that at no stage in applying $P_2^{-1}$ does one have to apply a representation of the inverse of the highly singular matrix $\mathcal{M}_{N_t}$, which is a key advantage of the preconditioner $P_2$ over $P_1$.

\vspace{0.5em}

We highlight that our methodology for constructing preconditioners of the form $P_1$ and $P_2$ may be readily tailored to the matrix system \eqref{for2system1} arising from the PDE arising from the alternative problem formulation discussed in \cref{sec::for2}.
\section{Discretisation using radial basis functions}
\label{sec::rbf}
The methods we have introduced so far are based on a finite difference discretisation of the partial differential equation. With practical imaging applications in mind, this means that the dimension of the discretised equation is typically proportional to the number of pixels in the image. Consequently, the number of degrees of freedom of the underlying equations is very large and can quickly become infeasible upon fine discretisation of the image. Assuming that the image of a now standard size for common smart phones is $3264\times 2448$ pixels, then solving an associated control problem with $100$ time-steps would lead to a Newton or Gauss--Newton system of dimensionality roughly $3$ billion.

\begin{figure}[htb!]
	\centering
	\subfloat[Initial image]
	{
		\includegraphics[height=4cm]{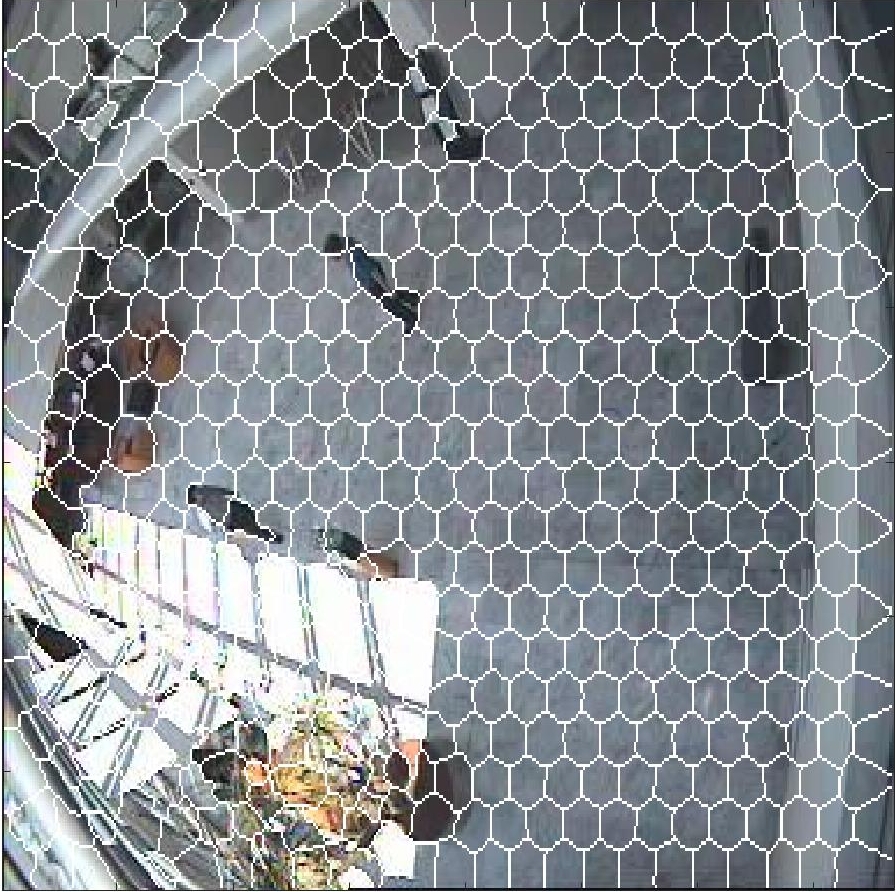} 
	}\hspace{1em}
	\subfloat[Target image]
	{
		\includegraphics[height=4cm]{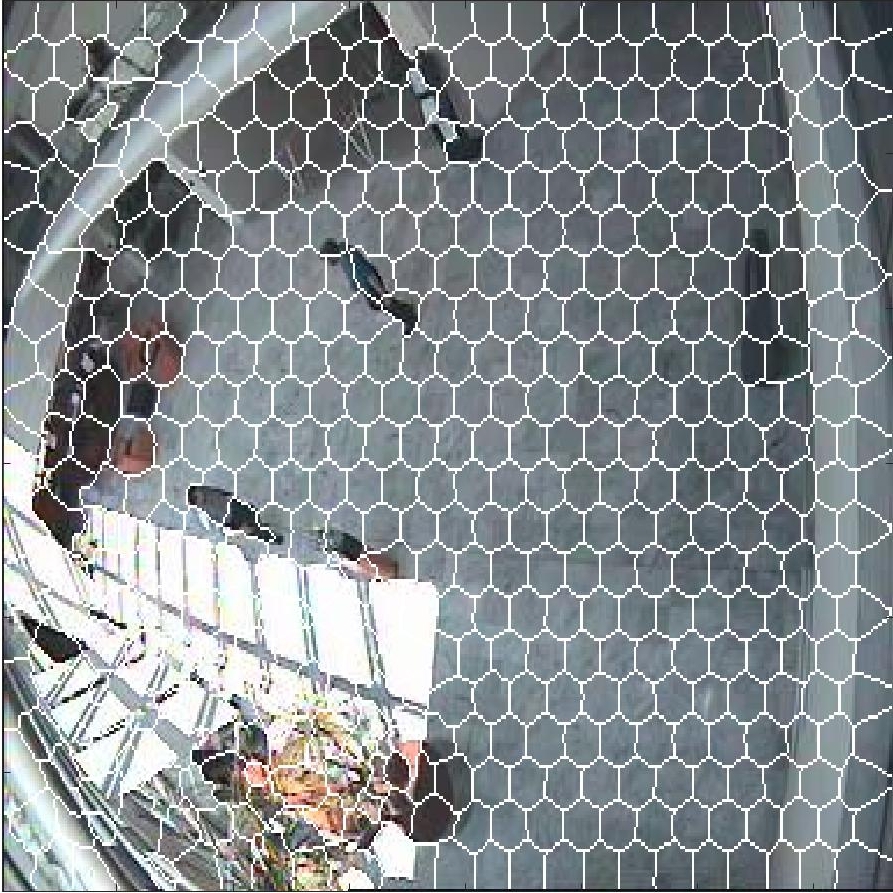} 
	}
	\caption{Image (left) used for superpixels with roughly $500$ superpixels, and image (right) with transferred superpixels. \label{fig::suppix}}
\end{figure}

In this section we wish to motivate an approach to reduce the number of degrees of freedom of the linear systems at the heart of the nonlinear iteration. Our technique is inspired by recent results on the use of reduction techniques based on clustered image information such as superpixels \cite{donne2015fast} or supervoxels \cite{amat2013fast}. Our aim is to reduce the complexity by applying a radial basis function approach, for which we create the scattered points as the centres of our superpixels, as illustrated in \cref{fig::suppix}\footnote{Two prototypical images taken from \url{http://homepages.inf.ed.ac.uk/rbf/CAVIARDATA1}.}. Before we discuss the detailed procedure for this discretisation approach, we point out that while the RBF methodology typically creates dense matrices, an image that is well described with a small number of superpixel will typically result in a small matrix representing the discretisation of the differential equation.

\subsection{RBF collocation for Newton system}

We now consider the optimise-then-discretise method discussed in \cref{sec::OTD}, and re-examine the Newton system obtained:
\begin{align*}
	\ (s_{y})_{_t}+\nabla\cdot (\vec{m} s_y)+\nabla\cdot (\vec{s}_{m}y)&=-\big(y_t+\nabla\cdot(\vec{m}y)\big),\\
	\ \beta \, Q^*Q\vec{s}_{m}-y\nabla s_{p}-s_y \nabla{}p&=-\big(\beta \, Q^* Q\vec{m}-y\nabla{}p\big),\\
	\ -(s_{p})_{_t}-\vec{s}_{m}\cdot\nabla{}p-\vec{m}\cdot\nabla s_{p}  &=-\big(-p_t-\vec{m}\cdot\nabla{}p\big).
\end{align*}
We wish to apply a meshfree method involving radial basis functions. The approach we use is \emph{straight collocation}, where the solution sought is a linear sum of RBFs multiplied by unknown coefficients, obtained by solving a matrix system. In more detail, at the $k$-th time-step we seek a solution for $s_y$, $s_{m_1}$, $s_{m_2}$, $s_p$ by substituting
\begin{align}
	\ (s_y)\big|_{t=t_k}=\sum_{j}Y_{j,k}\phi_j,\quad\quad(s_{m_1})\big|_{t=t_k}={}&\sum_{j}M_{j,k}^{(1)}\phi_j, 
	\notag
	\\
	\ (s_{m_2})\big|_{t=t_k}=\sum_{j}M_{j,k}^{(2)}\phi_j,\quad\quad\hspace{0.15em}(s_p)\big|_{t=t_k}={}&\sum_{j}P_{j,k}\phi_j,
	\label{eq:RBF_expansions}
\end{align}
into the Newton system. The coefficients $Y_{j,k}$, $M_{j,k}^{(1)}$, $M_{j,k}^{(2)}$, $P_{j,k}$ are unknowns, and $\phi_j$ denote the radial basis functions used. Each RBF is of the form $\phi=\phi(r)$ where $r=\big\|\vec{x}-\vec{\xi}\big\|$, with $\vec{x}$ the position vector and $\vec{\xi}$ the centre of the RBF. For our initial experiments, we use Gaussian functions $\phi(r)=e^{-cr^2}$ for a (positive) constant $c$ as our RBFs, though there are many other possibilities of functions with values at specified points solely depending on their distance from the centre \cite{Kansa,Wendland}.
Consider, for simplicity, the use of backward Euler for time derivatives. Then, at each time-step and RBF centre in space, the collocation procedure applied to the Newton system takes the following form:
\begin{align*}
	\ &\frac{1}{\tau}\left[\sum_{j}Y_{j,k}\phi_j-\sum_{j}Y_{j,k-1}\phi_j\right]+m_1\cdot\sum_{j}Y_{j,k}\frac{\partial\phi_j}{\partial{}x_1}+m_2\cdot\sum_{j}Y_{j,k}\frac{\partial\phi_j}{\partial{}x_2} \\
	\ &\quad\quad\quad\quad+\left(\frac{\partial{}m_1}{\partial{}x_1}+\frac{\partial{}m_2}{\partial{}x_2}\right)\cdot\sum_{j}Y_{j,k}\phi_j+\frac{\partial{}y}{\partial{}x_1}\cdot\sum_{j}M_{j,k}^{(1)}\phi_j+\frac{\partial{}y}{\partial{}x_2}\cdot\sum_{j}M_{j,k}^{(2)}\phi_j \\
	\ &\quad\quad\quad\quad+y\left(\sum_{j}M_{j,k}^{(1)}\frac{\partial\phi_j}{\partial{}x_1}+\sum_{j}M_{j,k}^{(2)}\frac{\partial\phi_j}{\partial{}x_2}\right)=-\big(y_t+\nabla\cdot(\vec{m}y)\big), \\
	\ &\beta{}Q^* Q\cdot\left[\begin{array}{cc}
			\sum_{j}M_{j,k}^{(1)}\phi_j \\ \sum_{j}M_{j,k}^{(2)}\phi_j \\
	\end{array}\right]-y\cdot\sum_{j}P_{j,k}\nabla\phi_j-\nabla{}p\cdot\sum_{j}Y_{j,k}\phi_j=-\big(\beta \, Q^* Q\vec{m}-y\nabla{}p\big), \\
	\ &-\frac{1}{\tau}\left[\sum_{j}P_{j,k+1}\phi_j-\sum_{j}P_{j,k}\phi_j\right]-\frac{\partial{}p}{\partial{}x_1}\cdot\sum_{j}M_{j,k}^{(1)}\phi_j-\frac{\partial{}p}{\partial{}x_2}\cdot\sum_{j}M_{j,k}^{(2)}\phi_j \\
	\ &\quad\quad\quad\quad-m_1\cdot\sum_{j}P_{j,k}\frac{\partial\phi_j}{\partial{}x_1}-m_2\cdot\sum_{j}P_{j,k}\frac{\partial\phi_j}{\partial{}x_2}=-\big(-p_t-\vec{m}\cdot\nabla{}p\big).
\end{align*}
The right hand sides in these equations are evaluated at the current iterate $(y,m_1,m_2,p)$, which is of the form \eqref{eq:RBF_expansions} as well.
The exception here is at the final time $T=1$, where there are additional terms are needed in order to take account of the term $\frac{1}{2\gamma}\int_{\Omega}(y(\vec{x},1)-y_1(\vec{x}))^2~\dOmega$ within the cost functional $\mathcal{E}(y,\vec{m})$, which we will include in our working.
Combining all the terms into a saddle point system of the form \eqref{SaddlePt} gives
\begin{equation}\label{RBFsystem}
	\ \left[\begin{array}{cccc}
			A_{y} & A_{y,m_1} & A_{y,m_2} & B_y^T \\
			A_{y,m_1} & A_{m_{11}} & A_{m_{12}} &  B_{m_1}^T \\
			A_{y,m_2} & A_{m_{21}} & A_{m_{22}} &  B_{m_2}^T \\
			C_y & C_{m_1} & C_{m_2} & 0 \\
	\end{array}\right]\left[\begin{array}{c}
			\mathbf{s}_{\mathbf{y}} \\
			\mathbf{s}_{\mathbf{m}_1} \\
			\mathbf{s}_{\mathbf{m}_2} \\
			\mathbf{s}_{\mathbf{p}} \\
	\end{array}\right]=\mathbf{b}.
\end{equation}
The matrices $C_y$ and $B_y$ take the form
\begin{equation*}
	\ C_y=\left[\begin{array}{cccc}
			C_y^{(1)} & & & \\
			-C_s & C_y^{(2)} & & \\
			& \ddots & \ddots & \\
			& & -C_s & C_y^{(N_t)} \\
	\end{array}\right],\quad\quad{}B_y=\left[\begin{array}{cccc}
			B_y^{(1)} & & & \\
			-C_s & B_y^{(2)} & & \\
			& \ddots & \ddots & \\
			& & -C_s & B_y^{(N_t)} \\
	\end{array}\right],
\end{equation*}
where the matrices $C_y^{(k)}$, $B_y^{(k)}$, and $C_s$ have entries
\begin{align*}
	\ \big(C_y^{(k)}\big)_{i,j}={}&\left(\phi_j+\tau m_1 \frac{\partial\phi_j}{\partial{}x_1}+\tau m_2 \frac{\partial\phi_j}{\partial{}x_2}+\tau\left(\frac{\partial{}m_1}{\partial{}x_1}+\frac{\partial{}m_2}{\partial{}x_2}\right)\phi_j\right)\Bigg|_{\vec{x}=\vec{x}_i}, \\
	\ \big(B_y^{(k)}\big)_{j,i}={}&\left(\phi_j-\tau m_1 \frac{\partial\phi_j}{\partial{}x_1}-\tau m_2 \frac{\partial\phi_j}{\partial{}x_2}\right)\Bigg|_{\vec{x}=\vec{x}_i}, \\
	\ \big(C_s\big)_{i,j}={}&\phi_j\big|_{\vec{x}=\vec{x}_i},
\end{align*}
with $m_1$, $m_2$ evaluated at the $k$-th time-step when $C_y^{(k)}$, $B_y^{(k)}$ are constructed. The remaining matrices $C_{m_1}$, $C_{m_2}$, $B_{m_1}$, $B_{m_2}$, $A_y$, $A_{m_1}$, $A_{m_2}$, $A_{y,m_1}$, $A_{y,m_2}$ are of block diagonal form:
\begin{align*}
	\ C_{m_1}={}&\text{blkdiag}\left(C_{m_1}^{(1)},...,C_{m_1}^{(N_t)}\right),\quad\quad\hspace{1.55em}{}C_{m_2}=\text{blkdiag}\left(C_{m_2}^{(1)},...,C_{m_2}^{(N_t)}\right), \\
	\ B_{m_1}={}&\text{blkdiag}\left(B_{m_1}^{(1)},...,B_{m_1}^{(N_t)}\right),\quad\quad\hspace{1.45em}{}B_{m_2}=\text{blkdiag}\left(B_{m_2}^{(1)},...,B_{m_2}^{(N_t)}\right), \\
	\ A_y={}&\text{blkdiag}\left(0,0,...,0,\gamma^{-1}C_s\right), \\
	\ A_{m_{11}}={}&\text{blkdiag}\left(A_{m_{11}}^{(1)},...,A_{m_{11}}^{(N_t)}\right),\quad\quad\hspace{0.8em}{}A_{m_{12}}=\text{blkdiag}\left(A_{m_{12}}^{(1)},...,A_{m_{12}}^{(N_t)}\right), \\
	\ A_{m_{21}}={}&\text{blkdiag}\left(A_{m_{21}}^{(1)},...,A_{m_{21}}^{(N_t)}\right),\quad\quad\hspace{0.8em}{}A_{m_{22}}=\text{blkdiag}\left(A_{m_{22}}^{(1)},...,A_{m_{22}}^{(N_t)}\right), \\
	\ A_{y,m_1}={}&\text{blkdiag}\left(A_{y,m_1}^{(1)},...,A_{y,m_1}^{(N_t)}\right),\quad\quad{}A_{y,m_2}=\text{blkdiag}\left(A_{y,m_2}^{(1)},...,A_{y,m_2}^{(N_t)}\right),% \\
%\ A_{m_1,y}={}&\text{blkdiag}\left(A_{m_1,y}^{(1)},...,A_{m_1,y}^{(N_t)}\right),\quad\quad{}A_{m_2,y}=\text{blkdiag}\left(A_{m_2,y}^{(1)},...,A_{m_2,y}^{(N_t)}\right),
\end{align*}
where
\begin{align*}
	\ \big(C_{m_1}^{(k)}\big)_{i,j}={}&\tau\left(\frac{\partial{}y}{\partial{}x_1}\phi_j+y\frac{\partial\phi_j}{\partial{}x_1}\right)\Bigg|_{\vec{x}=\vec{x}_i},\quad\quad\big(C_{m_2}^{(k)}\big)_{i,j}=\tau\left(\frac{\partial{}y}{\partial{}x_2}\phi_j+y\frac{\partial\phi_j}{\partial{}x_2}\right)\Bigg|_{\vec{x}=\vec{x}_i}, \\
	\ \big(B_{m_1}^{(k)}\big)_{j,i}={}&-\tau\left(y\frac{\partial\phi_j}{\partial{}x_1}\right)\Bigg|_{\vec{x}=\vec{x}_i},\quad\quad\hspace{2.875em}\big(B_{m_2}^{(k)}\big)_{j,i}=-\tau\left(y\frac{\partial\phi_j}{\partial{}x_2}\right)\Bigg|_{\vec{x}=\vec{x}_i}, \\
	\ \big(A_{m_{11}}^{(k)}\big)_{i,j}={}&\beta \, \tau\big((Q^* Q)_{1,1}\phi_j\big)\Big|_{\vec{x}=\vec{x}_i},\quad\quad\hspace{1.45em}\big(A_{m_{12}}^{(k)}\big)_{i,j}=\beta \, \tau\big((Q^* Q)_{1,2}\phi_j\big)\Big|_{\vec{x}=\vec{x}_i}, \\
	\ \big(A_{m_{21}}^{(k)}\big)_{i,j}={}&\beta \, \tau\big((Q^* Q)_{2,1}\phi_j\big)\Big|_{\vec{x}=\vec{x}_i},\quad\quad\hspace{1.45em}\big(A_{m_{22}}^{(k)}\big)_{i,j}=\beta \, \tau\big((Q^* Q)_{2,2}\phi_j\big)\Big|_{\vec{x}=\vec{x}_i}, \\
	\ \big(A_{y,m_1}^{(k)}\big)_{i,j}={}&-\tau\left(\frac{\partial{}p}{\partial{}x_1}\phi_j\right)\Bigg|_{\vec{x}=\vec{x}_i},\quad\quad\hspace{1.725em}\big(A_{y,m_2}^{(k)}\big)_{i,j}=-\tau\left(\frac{\partial{}p}{\partial{}x_2}\phi_j\right)\Bigg|_{\vec{x}=\vec{x}_i},
\end{align*}
with the relevant functions again evaluated at the $k$-th time-step. In each equation, the points $\vec{x}_i$ correspond to RBF centres chosen, and the vectors $\mathbf{s}_{\mathbf{y}}$, $\mathbf{s}_{\mathbf{m}_1}$, $\mathbf{s}_{\mathbf{m}_2}$, and $\mathbf{s}_{\mathbf{p}}$ concatenate the terms $Y_{j,k}$, $M_{j,k}^{(1)}$, $M_{j,k}^{(2)}$, and $P_{j,k}$, over all RBF centres and all time-steps. The terms $(Q^* Q)_{i,j}$, $i,j=1,2$, denote the $(i,j)$ blocks of the matrix $Q^* Q$. For the natural choices $Q=\text{blkdiag}(I,I)$ or $Q=\text{blkdiag}(\nabla,\nabla)$, the matrix $Q^* Q$ is given by $\text{blkdiag}(I,I)$ or $\text{blkdiag}(-\nabla^2,-\nabla^2)$, and in particular the matrices $A_{m_{12}}=A_{m_{21}}=0$. We also highlight that a Gauss--Newton approach for the discretised system, as discussed for the discretise-then-optimise method in \cref{sec::DTO}, would relate to the blocks $A_{y,m_1}$ and $A_{y,m_2}$ being the zero matrices.

One may consider preconditioners of the form $P_1$ and $P_2$, as derived in \cref{sec:Preconditioning}, for the system \eqref{RBFsystem}, provided one takes account of the (possibly dense) structure of the sub-blocks arising from the RBF collocation method.

\vspace{0.5em}

\textbf{Note.}~If the parameter $\delta>0$ within the cost functional $\mathcal{E}$, the only change arising in the above working would concern the matrix $A_y$, which would then be given by
\begin{equation*}
	\ A_y=\text{blkdiag}\left(\delta\tau{}C_s,\delta\tau{}C_s,...,\delta\tau{}C_s,(\gamma^{-1}+\delta\tau)C_s\right).
\end{equation*}

\subsection{Preconditioning}\label{sec:Preconditioning_RBF}

As the matrix systems that arise from the use of radial basis functions tend to be smaller and denser than those resulting from discretisation schemes such as finite differences and finite elements, there is in general less flexibility when designing fast and robust preconditioners. However, we believe that the construction of preconditioned iterative solvers remains useful, as one may therefore work with each time-step separately on a computer. This will decrease the storage requirements\deleted[id=RH]{ on a computer}, and will reduce the dimension of the matrices being solved for directly, as the sub-blocks arising from individual time-steps are relatively straightforward to solve for.

The idea for the preconditioner is the same as in \cref{sec:Preconditioning}, and we present an analogous preconditioner to $P_1$ as stated in \cref{sec::precond1}:
\begin{equation}
	\label{eq:P1_RBF}
	\ P_1=\left[\begin{array}{cccc}
			\widehat{A}_y & 0 & 0 & 0 \\
			0 & A_{m_{11}} & 0 & 0 \\
			0 & 0 & A_{m_{22}} & 0 \\
			0 & 0 & 0 & \widehat{S} \\
	\end{array}\right],
\end{equation}
where
\begin{align*}
	\ \widetilde{S}={}&C_y \widehat{A}_y^{-1} B_y^T + C_{m_1} A_{m_1}^{-1} B_{m_1}^T + C_{m_2} A_{m_2}^{-1} B_{m_2}^T \\
	\ \approx{}&\left(C_y+\mathcal{M}_1\right)\widehat{A}_y^{-1}\left(B_y^T+\mathcal{M}_2\right)=:\widehat{S}.
\end{align*}
Here, $\widehat{A}_y$ is an approximation of the often highly singular matrix $A_y$, $\widetilde{S}$ is generated by neglecting the off-diagonal terms in $(1,1)$ block, and $\mathcal{M}_1$, $\mathcal{M}_2$ are chosen such that
\begin{equation*}
	\ \mathcal{M}_1 \widehat{A}_y^{-1} \mathcal{M}_2 \approx C_{m_1} A_{m_1}^{-1} B_{m_1}^T + C_{m_2} A_{m_2}^{-1} B_{m_2}^T.
\end{equation*}
We can of course construct analogous block triangular preconditioners. To approximate the inverses of the sub-blocks of $\widehat{A}_y$, $A_{m_{11}}$, $A_{m_{22}}$, $C_y+\mathcal{M}_1$, $B_y^T+\mathcal{M}_2$ within the preconditioner, it is reasonable to apply direct solvers as the matrices are typically relatively small and dense. It is of course possible to replace this with a multigrid, domain decomposition, or other iterative scheme, as an approximation of the constituent matrices.

\section{Numerical results}\label{sec::Results}
We now present the results of a number of numerical experiments, making use both of the finite difference discretisation outlined in \cref{sec::discr}, and the radial basis function scheme of \cref{sec::rbf}.
All experiments are implemented in Matlab.

% \added[id=RH]{\textbf{It appears that the image data is not in the Git repository.}}
% \added[id=MS]{\textbf{Image data are now submitted, sorry.}}

\subsection{Finite difference discretisation}
\label{sec:FD}

\begin{figure}[htb!]
	\setlength\figureheight{0.2\linewidth} 
	\setlength\figurewidth{0.2\linewidth}
	\centering
	\subfloat[Linear solver iterations]{
		\includegraphics[height=4cm]{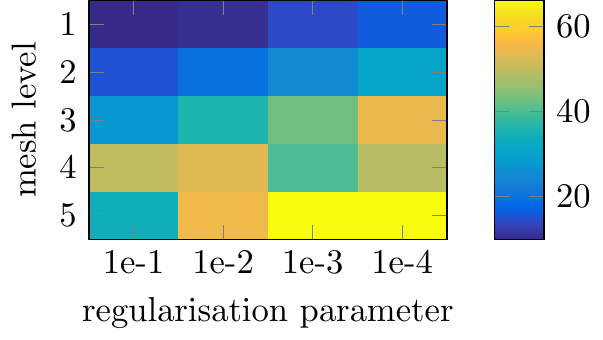} 
	} 
	\subfloat[Gauss--Newton iterations]{
		\includegraphics[height=4cm]{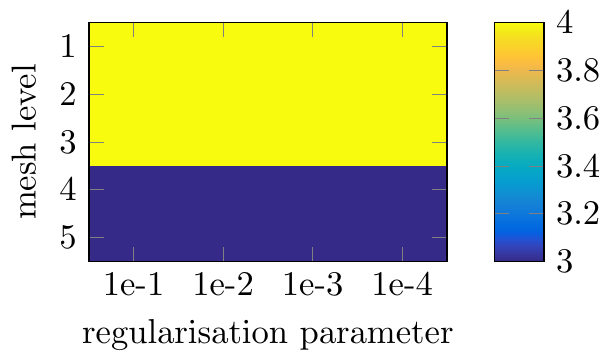} 
	}
	\caption{Results for the preconditioned iteration with $5$~mesh levels ranging from $n_x=2^{3}$ to $n_x=2^7$ degrees of freedom per spatial dimension. So the total number of spatial degrees of freedom ranges from $2^6$ to $2^{14}$. This variation is shown in the vertical axis, and the horizontal axis shows value of the regularisation parameter ranging from $10^{-1}$ to $10^{-4}.$ For the linear solver, in this case preconditioned CGS, we show iterations per Gauss--Newton step.}
	\label{fig:data_comp1}
\end{figure}

In our first test we consider problem \eqref{TranspCOST}--\eqref{TranspEQN} and employ the Gauss-Newton scheme \eqref{for1system1} in the finite difference setting. The parameters are chosen to be $\gamma=1$ and $Q=I.$ The regularisation parameter $\beta$ is typically varied in our experiments. We compare the performance of the preconditioner presented in \cref{sec::precond2} and recall that the $(1,1)$-block of the system matrix governing \eqref{for1system1} is highly singular. We use the same discretisation level in time as we use for the spatial domain, and compare $5$ spatial mesh levels for the synthetic data. For the image data we use ten time-steps and the same number of intermediate images. The image data is as depicted in \cref{fig::motivation}. We here use $100\times100$ pixel black and white images, where the values are scaled to be between zero and one. We choose the \cgs\ method \cite[Ch.~7.4.1]{book::Saad} as the iterative scheme for the Gauss--Newton system. This method is 
stopped when a certain tolerance (we use $10^{-6}$) for the relative residual norm of the linear system \eqref{for1system1} is reached, starting from an all-zero initial guess. Notice that in Matlab's implementation of \cgs, the residual is measured in the Euclidean norm. The outer Gauss--Newton scheme is stopped once the relative Euclidean distance between consecutive iterates falls below $10^{-4}$. \cref{fig:data_comp1} illustrates the average number of \cgs\ iterations on the one hand, and the number of Gauss--Newton steps on the other. The results are obtained using the preconditioner $P_2$, see \eqref{eq:prec1}, presented in \cref{sec::precond2}, and illustrate that this technique performs robustly with respect to the number of degrees of freedom and changes in the regularisation parameter.
%
% {In our first test we compare the performance of the preconditioner presented in \cref{sec::precond2} for the highly singular $(1,1)$-block. We compare $5$ mesh levels and choose the \cgs\ method {book::Saad} as the iterative scheme. The method is stopped when a tolerance of $10^{-6}$ for the linear system solver and a tolerance of $10^{-4}$ for the Gauss--Newton scheme are achieved. We measure the relative distance between new and old solution of the Gauss--Newton method. \cref{fig:data_comp1} illustrates the average number of \cgs\ iterations on the one hand, and the number of Gauss--Newton steps on the other. The results are obtained using the preconditioner $P_2$, and illustrate that this technique performs robustly with respect to the number of degrees of freedom and changes in the regularisation parameter.}

\begin{figure}[htb!]
	\centering
	\subfloat[Computed state at time step $5$]{
		\includegraphics[height=6cm]{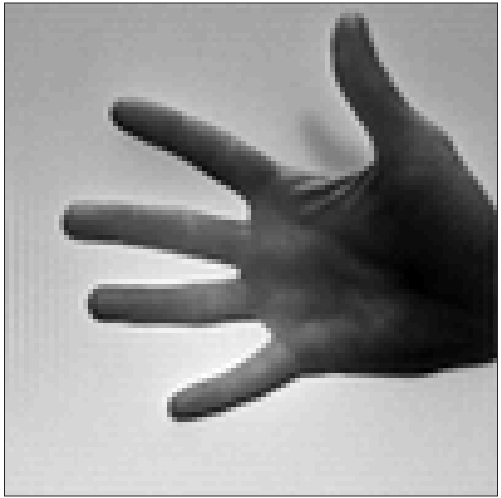} 
	} 
% \subfloat[Gauss--Newton Iterations]{
% \includegraphics[height=3cm]{Plots/hand_state6.png} 
% }
	\subfloat[Computed control at time step $5$]{
		\includegraphics[height=6cm]{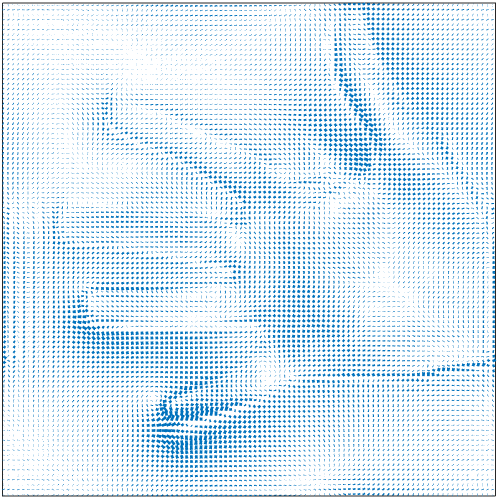} 
	}
	\setlength\figureheight{0.2\linewidth} 
	\setlength\figurewidth{0.5\linewidth}
	\subfloat[Varying $\delta$]{
		\includegraphics[height=4cm]{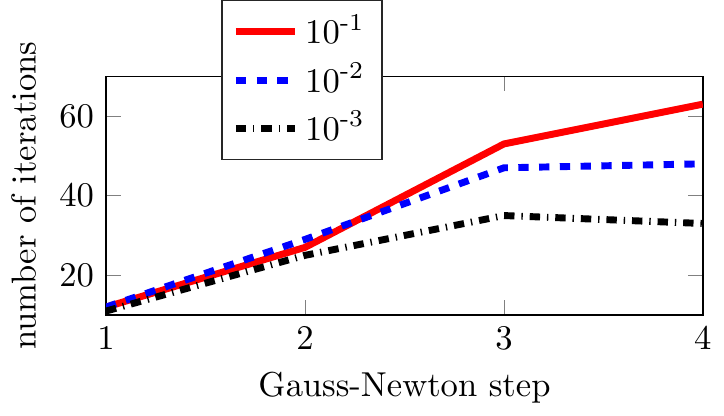} 		
	}
	\caption{Results for the preconditioned iteration for the optical flow problem. We show an instance of the control and the state as well as the number of CGS iterations per Gauss--Newton step for varying regularisation parameter $\delta$.}
	\label{fig:data_comp2}
\end{figure}

We also report results for the optical flow problem \cite{BIK02,borzi2003optimal}, i.e., we take the objective function to be given by \eqref{TranspOptical} with $\delta =10^{-3},$ $\beta=10^{-2},$ and $\gamma=10.$ The tolerances are set to be $10^{-3}$ and $10^{-2}$ for the linear and nonlinear solver, respectively.
%\begin{align}\label{TranspOptical}
%\begin{split}
%\E(y,\vec{m})={}&\frac{1}{2\gamma}\int_{\Omega}(y(\vec{x},1)-y_1(\vec{x}))^2~\dOmega+\frac{\delta}{2}\int_{0}^{1}\int_{\Omega}(y(\vec{x},t)-\bar{y}(\vec{x},t))^2\\
%&\quad\quad+\frac{\beta}{2}\int_{0}^{1}\int_{\Omega}(Q\vec{m}(\vec{x},t))^{2}~\dOmega \, \dt.
%\end{split}
%\end{align}
We here assume that $\bar{y}(\vec{x},1)=y_1$, and when discretised in time $\bar{y}(\vec{x},t_i)=y_i$ corresponds to a given image. As intermediate values for the desired state we chose the intermediate images from \cite{black1998eigentracking}. It is clear that this setup is covered by our previous discussion and the matrix representing the state contributions of the objection function $\M_{1,1}=\gamma^{-1}\M_{N_t}+\delta\tau\bar{\M}.$ One can readily apply the preconditioning techniques introduced in \cref{sec:Preconditioning}, and we consider the implementation of preconditioner $P_1$; see \eqref{eq:P1}. We show the results for our methodology in \cref{fig:data_comp2}. The desired state for this case was chosen as the hand sequence used in \cite{black1998eigentracking}.

We observe robustness with respect to the matrix dimension and the parameters involved in the problem setup, in both the Gauss--Newton iterations required and the number of steps of the preconditioned iterative method. This indicates the effectiveness of our preconditioning strategy.

\subsection{RBF technique}
\label{sec:RBF}

% \added[id=RH]{\textbf{Why do we apply it to artificial data and not to the same images as in \cref{sec:FD}??}}

To provide an indication of the applicability of radial basis function methods to the problems under consideration, we now provide details of the results obtained using this strategy.
We wish to test our proposed preconditioned iterative method on the matrix systems arising from our proposed RBF technique. For our next test problem, we provide a proof of concept that the methodology can be applied to dense matrix systems arising from Gaussian basis functions. We therefore select RBF centres using the superpixels computed using the initial image from \cref{fig::suppix}, and take the desired state to be a greyscale linear mapping, over the time variable, between this and the target image from the same figure. Within the cost functional \eqref{TranspOptical}, we take $\gamma=1$ and $Q=\text{blkdiag}(\nabla,\nabla)$. We test a range of values of $\delta$ and $\beta$, as well as time-steps $\tau=0.1$ and $\tau=0.05$ within the time interval $[0,1]$. In each case we apply outer (Gauss--Newton) iteration to a tolerance of $10^{-2}$ for the relative distance between new and old solution, and use a \gmres\ method preconditioned by $P_1$ (as described in \cref{sec:Preconditioning_RBF}; see \eqref{eq:P1_RBF}) to solve the matrix systems obtained from the outer iteration. Our Krylov solver is run to a tolerance of $10^{-4}$. As before, all norms are Euclidean in Matlab's implementation. The lowest number of outer iterations required for convergence was three, so to measure the effectiveness of the \gmres\ solver we present the average number of iterations for the first three Newton steps (as in general we find that the first Newton steps lead to the largest \gmres\ iteration counts). We present our results in \cref{RBFTable}.

We find that the \gmres\ iteration numbers are reasonably robust with respect to $\delta$ and $\tau$, but with greater dependence as $\beta$ than the finite difference approach. However, due to the fact that the matrix systems are much smaller for the radial basis function approach than for finite differences, the increase in achievable accuracy can compensate for the larger iteration counts. It would also be possible to apply our strategies using compactly supported RBFs, for instance Wendland functions \cite{Wendland}, instead of Gaussians, and thus exploit the sparsity of the resulting matrix systems.

\begin{table}[h]
	\renewcommand{\arraystretch}{1.1}
	\begin{footnotesize}
		\caption{Results for the radial basis function test problem, for $Q=\text{blkdiag}(\nabla,\nabla)$, $T=1$ and $\gamma=1$, using 10 time-steps (top) and 20 time-steps (bottom). Stated are the average number of \gmres\ iterations for the first three outer iterations, for a range of values of $\delta$ and $\beta$.}
		\label{RBFTable}
%\vspace{-0.5em}
		\begin{center}
			\begin{tabular}{|c||c|c|c|c|}
				\hline
				$\tau=0.1$ & $\beta=1$ & $\beta=10^{-1}$ & $\beta=10^{-2}$ & $\beta=10^{-3}$ \\ \hline \hline
				$\delta=0$ & 13 & 22.7 & 34.3 & 51.3 \\
				$\delta=0.1$ & 12.3 & 20.3 & 32.3 & 48.7 \\
				$\delta=1$ & 12.3 & 20 & 32.7 & 59.3 \\
				$\delta=10$ & 13 & 21.7 & 38.3 & 79.3 \\ \hline
			\end{tabular}

			\vspace{0.5em}

			\begin{tabular}{|c||c|c|c|c|}
				\hline
				$\tau=0.05$ & $\beta=1$ & $\beta=10^{-1}$ & $\beta=10^{-2}$ & $\beta=10^{-3}$ \\ \hline \hline
				$\delta=0$ & 16 & 25.7 & 51.7 & 60.3 \\
				$\delta=0.1$ & 14 & 23 & 52 & 62.3 \\
				$\delta=1$ & 14.3 & 22.3 & 42 & 71.3 \\
				$\delta=10$ & 14 & 24.7 & 43.3 & 74.3 \\ \hline
			\end{tabular}
		\end{center}
	\end{footnotesize}
\end{table}

\section{Conclusion}

In this paper, we have presented numerical methods for the solution of optimal transport problems arising in image metamorphosis. We have discussed the application of Newton and Gauss--Newton methods, using finite difference schemes and meshless methods for the discretisation of the optimality conditions. We presented fast and effective preconditioners which may be applied within Krylov subspace methods to solve the resulting matrix systems, with a focus on the large dimensions of the matrices when many time-steps are taken to solve the systems of PDEs. Encouraging numerical results indicate the potency of the solvers presented.

% 
% \bibliographystyle{siam}
% \bibliography{data}%JWP directory

\end{document}